\documentclass[11pt, final]{article} 

\usepackage{amsfonts,amsthm,amsmath,amssymb, mathtools}
\usepackage{dsfont}
\usepackage{graphicx,color}
\usepackage{hyperref}
\usepackage{a4wide}
\usepackage{multirow}
\usepackage{cleveref}
\usepackage{enumitem}
\usepackage{tikz-cd}
\usepackage{subcaption}

\newtheorem{proposition}{Proposition}[section]
\newtheorem{theorem}[proposition]{Theorem}
\newtheorem{lemma}[proposition]{Lemma}

\newtheorem{definition}[proposition]{Definition}
\newtheorem{remark}[proposition]{Remark}

\renewenvironment{proof}{\medskip\noindent{\textbf{Proof.}}%
  \hspace{1pt}}{\hspace{-5pt}{\nobreak\quad\nobreak\hfill\nobreak%
    $\square$\vspace{2pt}\par}\smallskip\goodbreak}

\newenvironment{proofof}[1]{\smallskip\noindent{\textbf{Proof~of~#1.}}%
  \hspace{1pt}}{\hspace{-5pt}{\nobreak\quad\nobreak\hfill\nobreak%
    $\square$\vspace{2pt}\par}\smallskip\goodbreak}

\numberwithin{equation}{section}
\numberwithin{figure}{section}
\numberwithin{table}{section}

\allowdisplaybreaks[4]

\usepackage{marginnote}

\renewcommand{\phi}{\varphi}
\renewcommand{\theta}{\vartheta}
\renewcommand{\epsilon}{\varepsilon}
\renewcommand{\L}[1]{\mathbf{L^#1}}

\renewcommand{\d}[1]{\mathinner{\mathrm{d}{#1}}}

\renewcommand{\L}[1]{\mathbf{L^#1}}
\newcommand{\dd}[1]{\mathinner{\mathrm{d}{#1}}}

\newcommand{\C}[1]{\mathbf{C^{#1}}}

\newcommand{\W}[2]{\mathbf{W^{#1,#2}}}
\newcommand{\BV}{\mathbf{BV}}

\newcommand{\modulo}[1]{{\left|#1\right|}}
\newcommand{\norma}[1]{{\left\|#1\right\|}}
\newcommand{\reali}{{\mathbb{R}}}

\newcommand{\interi}{{\mathbb{Z}}}
\newcommand{\Lip}{\mathop\mathbf{Lip}}
\newcommand{\tv}{\mathop\mathrm{TV}}
\newcommand{\essinf}{\mathop\mathrm{ess~inf}}
\newcommand{\esssup}{\mathop\mathrm{ess~sup}}

\newcommand{\caratt}[1]{\chi_{\strut#1}}

\newcommand{\rev}[1]{{\color{black} #1}}

\begin{document}

\title{General Stability Estimates in \\ NonLocal Traffic Models for
  Several Populations}

\author{R.M.~Colombo$^1$ \and M.~Garavello$^2$ \and C.~Nocita$^2$}

\maketitle

\footnotetext[1]{Unit\`a INdAM \& Dipartimento di Ingegneria
  dell'Informazione, Universit\`a di Brescia, Italy.\hfill\\
  \texttt{rinaldo.colombo@unibs.it}}

\footnotetext[2]{Dipartimento di Matematica e Applicazioni,
  Universit\`a di Milano--Bicocca, Italy.\hfill\\
  \texttt{mauro.garavello@unimib.it} and
  \texttt{c.nocita@campus.unimib.it}}

\begin{abstract}
  \noindent We prove global existence, uniqueness and $\L1$ stability
  of solutions to general systems of nonlocal conservation laws
  modeling multiclass vehicular traffic. Each class follows its own
  speed law and has specific effects on the other classes'
  speeds. Moreover, general explicit dependencies of the speed laws on
  space and time are allowed. Solutions are proved to depend
  continuously --- in suitable norms --- on all terms appearing in the
  equations, as well as on the initial data. Numerical simulations
  show the relevance and the effects of the nonlocal terms.

  \medskip


  \medskip

  \noindent\textbf{MSC~2020:} 35L65, 76A30, 65M22
\end{abstract}

\section{Introduction}
\label{sec:introduction}

In this paper we analyze a nonlocal multiclass, or multipopulation,
macroscopic traffic flow model. Here, $\rho_1, \ldots,\rho_n$ are the
densities of $n$ different populations (or classes) of drivers or
vehicles, whose evolution is described by the following system of
nonlinear conservation laws:
\begin{equation}
  \label{eq:9}
  \left\{
    \begin{array}{ll}
      \partial_t \rho_i
      +
      \partial_x \left(\rho_i\; v_i(t, x, \eta_i * \rho)\right)=0
      & (t,x) \in [t_o, +\infty\mathclose[ \times \reali
      \\
      \rho_i(t_o,x)=\rho_i^o(x)
      & x \in \reali
    \end{array}
  \right.
  \qquad
  i = 1, \ldots, n\,,
\end{equation}
where $\eta_1, \ldots, \eta_n \colon \reali \rightarrow \reali^n$ are
suitable weights,
$v_1, \ldots, v_n \colon \reali_+ \times\reali \times \reali^n
\rightarrow \reali$ are the time and space dependent speed laws, while
$\rho_1^o, \ldots, \rho_n^o$ are the initial densities,
\rev{and $\rho = (\rho_1, \ldots, \rho_n)$}. The present
setting allows, in a time dependent framework, to account for space
inhomogeneities as well as for the possible different natures of
drivers/vehicles and of their interactions. Indeed, the equations
in~\eqref{eq:9} are coupled through the nonlocal terms $\eta_i * \rho$
where each of the space convolutions
\begin{displaymath}
  (\eta_i * \rho)_j \, (t,x)
  \coloneqq
  \left(\eta_{ij} * \rho_j (t)\right)(x)
  =
  \int_{\reali} \eta_{ij} (\xi) \; \rho_j (t,x-\xi) \d\xi
  \qquad
  i,j = 1, \ldots, n
\end{displaymath}
describes how the $i$-th population interacts with the $j$-th one.
Note that we consider the possibility of vehicles possessing a
non-standard visual horizon, as is the case with autonomous vehicles,
which typically have a comprehensive understanding of the traffic environment,
both ahead and behind.

Within this general framework, defined only by the two rather simple
assumptions~\ref{item:1} and~\ref{item:2},
\rev{defined in Section~\ref{sec:analytical-results},} besides proving global in
time existence and uniqueness of the solution to~\eqref{eq:9}, we show
the $\L1$-Lipschitz continuous dependence of solutions on the initial
data, their $\L1$ local Lipschitz continuity in time and also provide
detailed stability estimates on the dependence of solutions on $v$ and
on $\eta$\rev{; see also \cite{zbMATH07030125}}. Moreover, \emph{a priori} estimates on the total variation
in space of the solution are obtained.

The current literature offers several related results, due to the
interest in nonlocal traffic models. Well-posedness and stability of
multiclass systems with a coupling in the nonlocal terms were
considered in~\cite{MR2921871, MR3057143}. In particular, the latter
work is set mainly in a Radon measure setting with stability estimates
expressed by means of Wasserstein distance, thus requiring initial
data to have the same total mass.

Motivated by the effects of moving obstacles on vehicular traffic,
\cite{MR4442432} considers a coupled ODE-PDE system. Well-posedness of
weak solutions is shown for sufficiently small times and, under
stronger conditions on the convolution kernels, long time existence is
also proved.


The development of \emph{ad hoc} numerical methods for equations of
the type~\eqref{eq:9} was considered, for instance,
in~\cite{MR3342191, MR3759879, MR4362532}. Recall that nonlocal
systems fitting in the form~\eqref{eq:9} can be motivated also by
entirely different physical settings: see for
instance~\cite{MR4074000} devoted to the dynamics of bolts on a
conveyor belt, \cite{MR3397337} devoted to the flow of melted metal
and~\cite{Armbruster2006933, MR2679644} modeling supply chains. Crowd
dynamics is a further widely considered application, see for
instance~\cite{MR3460619, MR2902155}. A relaxation representation of
nonlocal conservation laws was considered in~\cite{MR4110434}.

Nonlocal conservation laws have been at the center of various papers
dealing with the nonlocal to local limit, e.g.~\cite{MR3342191,
  MR4651679, zbMATH07658225, MR3961295, MR4613802, MR3944408}.
In particular, the
negative result obtained in~\cite{MR3961295} is consistent with the
stability estimate on the dependence of solution on $\eta$ obtained
hereafter, which requires a rather strong norm.

Particular features of the present analytical setting are: the full
nonlinearity of~\eqref{eq:9}, possible thanks to the generality of the
functions $v_i$, the explicit dependence of the speed laws on the
space variable $x$ and on time $t$, the variety of the interactions
between the different populations allowed by the, possibly different,
$n^2$ functions $\eta_{ij}$.

Due to the nonlinearity of~\eqref{eq:9}, a natural tool is the Banach
Fixed Point Theorem, which we apply in
$\C0 \left([0,T]; \L1 (\reali; \reali^n)\right)$ and relies on the
representation of solutions to renewal equations by means of
characteristics, often referred to also as Lagrangian solutions, as
in~\cite{MR3057143}. With this tool, we prove existence, uniqueness
and $\L1$-Lipschitz continuous dependence of the solution on the
initial datum. Careful $\BV$ bounds allow first to extend the solution
globally in time. Then, we also get the local in time $\L1$-Lipschitz
continuous dependence of the solution on $t$, $v$ and $\eta$.

We stress that the stability estimates proved below require that the total
variation in space of the initial datum, and hence of the solution, be
finite. In the case of the (local) Lipschitz continuity in time this
is shown by the following elementary example:
\begin{displaymath}
  \left\{
    \begin{array}{l}
    \partial_t \rho + \partial_x \rho = 0
      \\
    \rho (0,x) = \rho_o^m (x)
    \end{array}
  \right.
  \quad \mbox{ where } \quad
  \rho_o^m (x) =
  \sum_{i=1}^m \caratt{[\frac{2i}{2m},\frac{2i+1}{2m}]} (x)
\end{displaymath}
which yields, for $\epsilon \in \mathopen]0, 1/ (2m)\mathclose[$,
\begin{displaymath}
  \norma{\rho (t+\epsilon) - \rho (t)}_{\L1 (\reali; \reali)}
  = \tv (\rho_o^m) \; \epsilon \,.
\end{displaymath}
In turn, the estimates on the total variation in space depend on the
$\W2\infty$ norm of the convolution kernel. On the one hand, this
motivates our assumptions being more restrictive than those
in~\cite{MR3057143}, where only continuity in time (not Lipschitz
continuity) is proved. On the other hand, this is also consistent with
the impossibility, in general, of letting $\eta$ (formally) converge
to a Dirac delta, which is equivalent to pass from a nonlocal to a
local problem, see~\cite{MR4651679, MR3961295, MR3944408}. We refer for instance to~\cite{MR3959426} for $\L{p}$ continuity estimates, with $p>1$.

Various other papers use fixed point arguments together with implicit
or explicit expressions for the solutions, see~\cite{MR2679644,
  MR3057143, MR3670045, MR3944408}. Alternatively, sequences of
approximate solutions constructed by means of numerical algorithms can
also be used to prove the existence of solutions to systems fitting
into~\eqref{eq:9}, see for instance~\cite{aggarwal2023systems,
  aggarwal2023wellposedness, MR3808157, MR3959349}.

Below, numerical integrations show specific features of the
model~\eqref{eq:9}, such as the effects of different horizons and
overtakes. Moreover, we compare a solution to the nonlocal
equation~\eqref{eq:9} with that of the classical
Lighthill-Whitham~\cite{LighthillWhitham} and Richards~\cite{Richards}
model in presence of space inhomogeneities (bottleneck).

The paper is organized as follows. In
Section~\ref{sec:analytical-results} we present the main analytical
results. Numerical integrations of~\eqref{eq:9} are shown in
Section~\ref{sec:numer-integr}. All proofs are deferred to
Section~\ref{sec:technical-details}.

\section{Analytical Results}
\label{sec:analytical-results}

Below we denote by $\partial_t$, $\partial_x$ and $\nabla_\rho$ the
partial derivatives with respect to, respectively, $t$, $x$ and
$(\rho_1, \ldots, \rho_n)$. \rev{Moreover, if $v = v (x,\rho)$ with
  $x \in \reali$ and $\rho \in \reali^n$, then the notation
  $[\partial_x v \quad \nabla_\rho v]$ is the vector
  $(\partial_x v, \partial_{\rho_1} v, \ldots, \partial_{\rho_n} v)$.}
Fix $t_o \in \rev{\reali_+ = [0, + \infty\mathclose[}$. We
study~\eqref{eq:9} under the following assumptions:
\begin{enumerate}[label=\textbf{(v)}, ref=\textup{\textbf{(v)}}]
\item \label{item:1}
  $v \, \rev{ = (v_1, \ldots, v_n)} \in \C0 ([t_o,
  +\infty\mathclose[;\mathcal{V}^n)$, where
  \begin{displaymath}
    \mathcal{V}
    \coloneqq
    \left\{v \colon \reali\times \reali^n \to \reali \colon
      v (\cdot, 0) \in \L\infty (\reali;\reali),\; [
      \partial_x v \quad \nabla_\rho v]
      \in \W2\infty (\reali\times\reali^n; \reali\times\reali^n)
    \right\} \,.
  \end{displaymath}
\end{enumerate}
\begin{enumerate}[label=\textbf{($\mathbf{\eta}$)},
  ref=\textup{($\mathbf{\eta}$)}]
\item \label{item:2} For $i=1,\ldots,n$,
  $\eta_i \in \W2\infty (\reali;\reali^n)$.
\end{enumerate}
\noindent Referring to $\mathcal{V}$ and $\mathcal{V}^n$, we use the
norms
\begin{displaymath}
  \begin{array}{c}
    \begin{array}{rcl}
      \norma{v_i}_{\mathcal{V}}
      & \coloneqq
      & \norma{v_i (\cdot, 0)}_{\L\infty (\reali; \reali)}
        +
        \norma{[\partial_x v_i \quad \nabla_\rho v_i]}_{\W2\infty (\reali\times\reali^n; \reali\times\reali^n)}
      \\[6pt]
      \norma{v}_{\mathcal{V}^n}
      & \coloneqq
      & \sqrt{\sum_{i=1}^n  \norma{v_i}_{\mathcal{V}}^2}
      \\[6pt]
      \norma{v}_{\C0([t_o,t];\mathcal{V}^n)}
      & \coloneqq
      & \sup_{\tau \in
        [t_o,t]}\norma{v (\tau)}_{\mathcal{V}^n}\,.
    \end{array}
  \end{array}
\end{displaymath}

\begin{remark}
  \label{rem:no3x}
  {\rm As the proofs below show, the third derivative
    $\partial^3_{xxx }v$ is not even required to exist. Nevertheless,
    we state assumption~\ref{item:1} and define
    $\norma{\,\cdot\,}_{\mathcal{V}}$ as above merely for simplicity.}
\end{remark}

As a first step, we formalize what we mean by solution
to~\eqref{eq:9}.

\begin{definition}
  \label{def:solution}
  Fix a non empty real interval $I$ with $\min I = t_o$. By
  \emph{solution to~\eqref{eq:9} on the time interval $I$} we mean a
  map
  $\rho = (\rho_1, \ldots, \rho_n) \in \C0 (I; \L1 (\reali;
  \reali^n))$ such that, for $i=1, \ldots, n$, setting
  \begin{equation}
    \label{eq:14}
    w_i (t,x) \coloneqq v_i\left(t, x,
      \left(\eta_{i1} * \rho_1 (t)\right) (x),
      \left(\eta_{i2} * \rho_2 (t)\right) (x),
      \ldots,
      \left(\eta_{in} * \rho_n (t)\right) (x)
    \right) \,,
  \end{equation}
  $\rho_i$ is a solution to
  \begin{equation}
    \label{eq:20}
    \left\{
      \begin{array}{l@{\qquad}l}
        \partial_t \rho_i + \partial_x \left(\rho_i \; w_i (t,x)\right) = 0
        & (t,x) \in I \times \reali
        \\
        \rho_i (t_o,x) = (\rho_o)_i (x)
        & x \in \reali \,.
      \end{array}
    \right.
  \end{equation}
\end{definition}

\noindent Above, by \emph{solution to~\eqref{eq:20}} we mean a
distributional solution~\cite[Definition~4.2]{MR1816648}, which is
also a weak entropy solution in the sense
of~\cite[Definition~1]{MR267257}, see~\cite[Lemma~5]{MR4371486},
\cite[Theorem~2.7]{MR3057143} or~\cite[Corollary~II.1]{MR1022305}.

We are now ready to state the main analytical result. To this aim, we
denote by $C (\cdots)$ a locally bounded quantity dependent on its
arguments, whose actual value is specified in the proofs.

\begin{theorem}
  \label{thm:main}
  Let~\ref{item:1} and~\ref{item:2} hold. Then, problem~\eqref{eq:9}
  generates a unique map
  \begin{displaymath}
    \mathcal{P}
    \colon
    [t_o, +\infty\mathclose[ \times
    \reali_+ \times
    (\L1 \cap \BV)(\reali; \reali^n)
    \to (\L1 \cap \BV)(\reali; \reali^n)
  \end{displaymath}
  with the following properties:
  \begin{enumerate}[label=\textup{\textbf{($\mathcal{P}$\arabic*)}},
    ref=\textup{\textbf{($\mathcal{P}$\arabic*)}}]
  \item \label{item:16} $\mathcal{P}$ is a process: for all
    $\bar t \geq t_o$ and $t',t'' \in \reali_+$,
    $\mathcal{P}_{\bar t,0}$ is the identity and
    $\mathcal{P}_{\bar t + t', t''} \circ \mathcal{P}_{\bar t, t'} =
    \mathcal{P}_{\bar t, t'+t''}$.
  \item \label{item:4} For any
    $\rho_o \in (\L1 \cap \BV)(\reali; \reali^n)$ and any $T>t_o$, the
    orbit $t \mapsto \mathcal{P}_{t_o,t} \rho_o$ is the unique global
    solution to~\eqref{eq:9} on $[t_o,\,T]$ in the sense of
    Definition~\ref{def:solution}.

  \item \label{item:5} $\mathcal{P}$ is locally Lipschitz continuous
    in the initial datum: for all $t \geq t_o$ and
    \rev{$\rho_o, \hat \rho_o \in (\L1\cap \BV) (\reali; \reali^n)$}
    \begin{eqnarray*}
      &
      & \norma{\mathcal{P}_{t_o,t}\rev{\rho_o}
        - \mathcal{P}_{t_o,t} \rev{\hat \rho_o}}_{\L1 (\reali; \reali^n)}
      \\
      & \leq
      & C\left(\norma{\eta}_{\W2\infty (\reali;\reali^{n\times n})},
        \norma{v}_{\C0([t_o,t];\mathcal{V}^n)},
        \norma{\rev{\rho_o}}_{\L1(\reali;\reali^n)},
        \norma{\rev{\hat\rho_o}}_{\L1 (\reali;\reali^)},
        \tv (\rev{\hat \rho_o}), t-t_o\right)
      \\
      && \times
         \norma{\rev{\rho_o - \hat \rho_o}}_{\L1 (\reali; \reali^n)} \,.
    \end{eqnarray*}
  \item \label{item:8} $\mathcal{P}$ is locally Lipschitz continuous
    in $t$: for all $T>t_o$, $t',t'' \in [t_o,T]$ and
    $\rho_o \in (\L1\cap \BV) (\reali; \reali^n)$
    \begin{eqnarray*}
      &
      & \norma{\mathcal{P}_{t_o,t'}\rho_o - \mathcal{P}_{t_o,t''} \rho_o}_{\L1 (\reali; \reali^n)}
      \\
      & \leq
      & C\left(\norma{\eta}_{\W2\infty (\reali;\reali^{n\times n})},
        \norma{v}_{\C0([t_o,T];\mathcal{V}^n)}, \norma{\rho_o}_{\L1 (\reali;\reali^n)}, \tv (\rho_o) , T-t_o\right)
        \modulo{t' {-} t''} \,.
    \end{eqnarray*}

  \item \label{item:6} Let \rev{$v, \hat v$} satisfy~\ref{item:1} and
    call \rev{$\mathcal{P},\hat {\mathcal{P}}$} the corresponding
    processes. Then, for all $t \geq t_o$ and
    $\rho_o \in (\L1\cap\BV) (\reali;\reali^n)$
    \begin{eqnarray*}
      &
      & \norma{\mathcal{P}_{t_o,t}\rho_o - \hat{\mathcal{P}}_{t_o,t}\rho_o}_{\L1 (\reali;\reali^n)}
      \\
      & \leq
      &  C\left(\norma{\eta}_{\W2\infty (\reali;\reali^{n\times n})},
        \norma{v}_{\C0([t_o,t];\mathcal{V}^n)},
        \norma{\rho_o}_{\L1(\reali;\reali^n)}, \tv(\rho_o) , t-t_o\right)
      \\
      &
      & \qquad \times
        \norma{v-\hat v}_{\C0([t_o,t];\mathcal{V}^n)} \; (t-t_o) \,.
    \end{eqnarray*}

  \item \label{item:7} Let \rev{$\eta, \hat \eta$}
    satisfy~\ref{item:2} and call
    \rev{$\mathcal{P}, \hat{\mathcal{P}}$} the corresponding
    processes. Then, for all $t \geq t_o$ and
    $\rho_o \in (\L1\cap\BV) (\reali;\reali^n)$
    \begin{eqnarray*}
      &
      & \norma{\mathcal{P}_{t_o,t}\rho_o - \hat{\mathcal{P}}_{t_o,t}\rho_o}_{\L1 (\reali;\reali^n)}
      \\
      & \leq
      &
        C\left({\norma{\eta}_{\W2\infty(\reali;\reali^{n \times n})} ,
        \norma{\hat \eta}_{\W2\infty(\reali;\reali^{n \times n})}},
        \norma{v}_{\C0([t_o,t];\mathcal{V}^n)},
        \norma{\rho_o}_{\L1(\reali;\reali^n)}, \tv(\rho_o), t-t_o\right)
      \\
      &
      & \qquad\times
        \norma{\eta-\hat \eta}_{\W2\infty(\reali;\reali^{n \times n})} \; (t-t_o) \,.
    \end{eqnarray*}

  \item \label{item:15} For all $t \geq t_o$ and
    $\rho_o \in (\L1\cap\BV) (\reali;\reali^n)$
    \begin{eqnarray*}
      \tv\left(\mathcal{P}_{t_o,t} \rho_o\right)
      & \leq
      & \left(
        \tv (\rho_o) +
        C \left(
        \norma{\eta}_{\W2\infty (\reali; \reali^{n\times n})},
        \norma{v}_{\C0([t_o,t];\mathcal{V}^n)}, \norma{\rho_o}_{\L1 (\reali;\reali^n)}
        \right) (t-t_o)
        \right)
      \\
      &
      & \quad \times
        \exp \left(C \left(
        \norma{\eta}_{\W2\infty (\reali; \reali^{n\times n})},
        \norma{v}_{\C0([t_o,t];\mathcal{V}^n)},
        \norma{\rho_o}_{\L1 (\reali;\reali^n)}
        \right) (t-t_o)\right) \,.
    \end{eqnarray*}

  \item \label{item:14} For every $i = 1, \ldots, n$, if
    $(\rho_o)_i \geq 0$, then for all $t \geq t_o$,
    $(\mathcal{P}_{t,t_o}\rho_o)_i \geq 0$.
  \end{enumerate}
\end{theorem}

\noindent The proof is deferred to
Section~\ref{sec:technical-details}. \rev{Since $\rho\equiv 0$ is a
  solution, the bounds~\ref{item:5} and~\ref{item:15} ensure also an
  $\L\infty$ bound on $\rho$.}

\section{Numerical Integrations}
\label{sec:numer-integr}

In the numerical integrations below, motivated by vehicular traffic,
we choose non negative values of each component of the initial density
$\rho_o$. Hence, thanks to~\ref{item:14}, it is sufficient that we
define for all $i=1,\ldots,n$ the speed laws $v_i$ only on
$\reali_+ \times \reali \times \reali_+^n$. Then, one can extend them
by regularity to all $\reali \times \reali^n$ in order to formally
recover the hypothesis~\ref{item:1}.

The kernel functions $\eta_{ij}$ describe how far the $i$-th
population \emph{``sees''} the $j$-th population. For simplicity, we
standardize
their choices as follows:\\
\begin{minipage}{0.4\linewidth}
    \includegraphics[width=\textwidth,
  keepaspectratio]{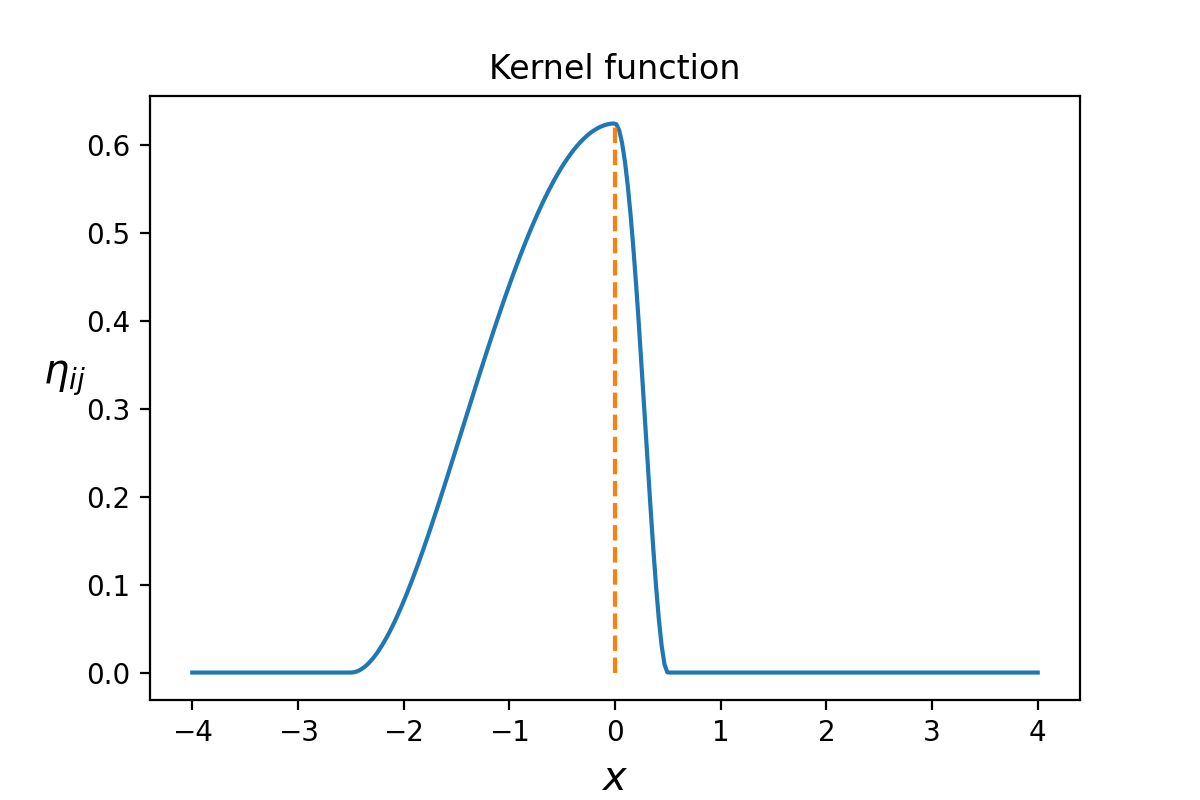}
\end{minipage}%
\begin{minipage}{0.6\linewidth}
\begin{equation}
  \label{eq:kernel}
  \eta_{ij} (x)
  \coloneqq
  \left\{
    \begin{array}{@{}l@{\qquad}r@{\,}c@{\,}l@{}}
      0
      & x
      & <
      & -f_{ij}
      \\
      A_{ij} \; \left(1-\left(\frac{x}{f_{ij}}\right)^2\right)^2
      & x
      & \in
      & [-f_{ij}, 0]
      \\
      A_{ij} \; \left(1-\left(\frac{x}{b_{ij}}\right)^2\right)^2
      & x
      & \in
      & [0,b_{ij}]
      \\
      0
      & x
      & >
      & b_{ij}
    \end{array}
  \right.
\end{equation}
\end{minipage}\\
where the normalization constants $A_{ij}$ are chosen so that
$\int_{\reali} \eta_{ij} (x) \d{x} = 1$. The forward, respectively
backward, horizons are $f_{ij}$, respectively $b_{ij}$, and are
typically chosen so that $f_{ij} \gg b_{ij}$, accounting for the far
greater relevance of the forward vision.

Note however that whenever a population of autonomous vehicles (AVs)
is present, it is reasonable to assume that they have information
about the whole traffic situation, also out of the standard visual
horizon, both in front and behind them. This is consistent with our
allowing the kernel function $\eta_{ij}$ to have, possibly, unbounded
support. At the same time, the backward horizon for a standard vehicle
may well be of the order of a vehicle length, coherently with the
little relevance of the behind traffic.

Throughout, we use the Lax-Friedrichs method~\cite[\S~4.6]{MR1925043}
adapted to deal with nonlocal fluxes. Indeed, fix a spatial mesh size
$\Delta x > 0$, a final time $T>0$, and construct the grid points
$x_k\coloneqq k \, \Delta x$, for $k \in \interi$, and, for
$m=1,\ldots, N_T$, the times $t_m\coloneqq t_{m-1} + \Delta t_m$,
where $t_0 \coloneqq 0$ and $\Delta t_m > 0$ to be specified below and
such that $T=\sum_{m=0}^{N_T}\Delta t_m$.  Call, for all
$k \in \interi$ and $m=0,\ldots, N_T-1$, the cell
$C_k\coloneqq\left[x_{k-\frac{1}{2}}, x_{k+\frac{1}{2}}\right]$, where
$x_{k+\frac{1}{2}}\coloneqq\left( x_{k+1}+x_k\right)/2$, and the time
interval $C^m\coloneqq[t_m, t_{m+1}\mathclose[$.  Then, for all
$i=1,\ldots,n$, the solution $\rho_i$ to~\eqref{eq:9} is approximated
by the piecewise constant function
$\sum_{k,m} (\rho_i)_k^m \, \caratt{C^m} (t) \; \caratt{C_k} (x)$
where, for $m=0,\ldots, N_T-1$ and $ k \in \interi$, $(\rho_i)_k^m$
is given by the Lax-Friedrichs scheme
\begin{equation*}
  \begin{cases}
    (\rho_i)_k^{m+1} = \frac{1}{2} \left( (\rho_i)_{k-1}^{m} + (\rho_i)_{k+1}^{m}\right) - \frac{\Delta t_m}{2\Delta x}\left[ F\left( (\rho_i)_{k+1}^{m} \right) - F\left((\rho_i)_{k-1}^{m}\right) \right]\\
    (\rho_i)_k^0= (\rho_o)_i(x_k),
  \end{cases}
\end{equation*}
with the numerical fluxes
$F((\rho_i)_k^m)\coloneqq (\rho_i)_k^m (v_i)_k^m$.  Here, the values
$(v_i)_k^m$ are defined by
\begin{eqnarray*}
  (v_i)_k^m
  & \coloneqq
& v_i\left(t_m, x_{k}, \Delta x \sum_{p \in \interi} \eta_{i1}(x_p)(\rho_1)^m_{k-p}, \ldots, \Delta x \sum_{p \in \interi} \eta_{in}(x_p)(\rho_n)^m_{k-p}\right),
\end{eqnarray*}
which are an approximation of
$v_i\left(t_m, x_{k}, \eta_{i1}*\rho_1(t_m,x_{k}), \ldots,
  \eta_{in}*\rho_n(t_m,x_k)\right)$.  For each
$m=0,\ldots, N_T-1$, $\Delta t_m $ is chosen to satisfy the CFL
condition~\cite[\S~4.4]{MR1925043}
\begin{equation}
  \label{CFL}
  \Delta t_m \leq \dfrac{\Delta x}{\displaystyle\max_{\substack{i=1,\ldots, n, \\ k\in \interi}} \{ (v_i)^m_k\}}.
\end{equation}
\rev{Finally, at the numerical level, we employed free flow boundary conditions.}

\subsection{The Role of the Horizon}
\label{subsec: horizon}

Here we consider $n=2$ populations of drivers whose unique difference
is their visual horizon. Such distinction is carried by the kernel
functions $\eta_{ij}$ in~\eqref{eq:kernel}.  More precisely, for
$i,j=1,2$, fix $b_{ij} = 0.01$, $f_{1j} = 1.5$ and $f_{2j} = 0.3$ so
that, forward, the first population is able to see farther than the
second one.

\begin{figure}[!ht]
  \begin{subfigure}{0.33\textwidth}
    \centering \includegraphics[width=1\textwidth, trim=20 10 20 10
    keepaspectratio]{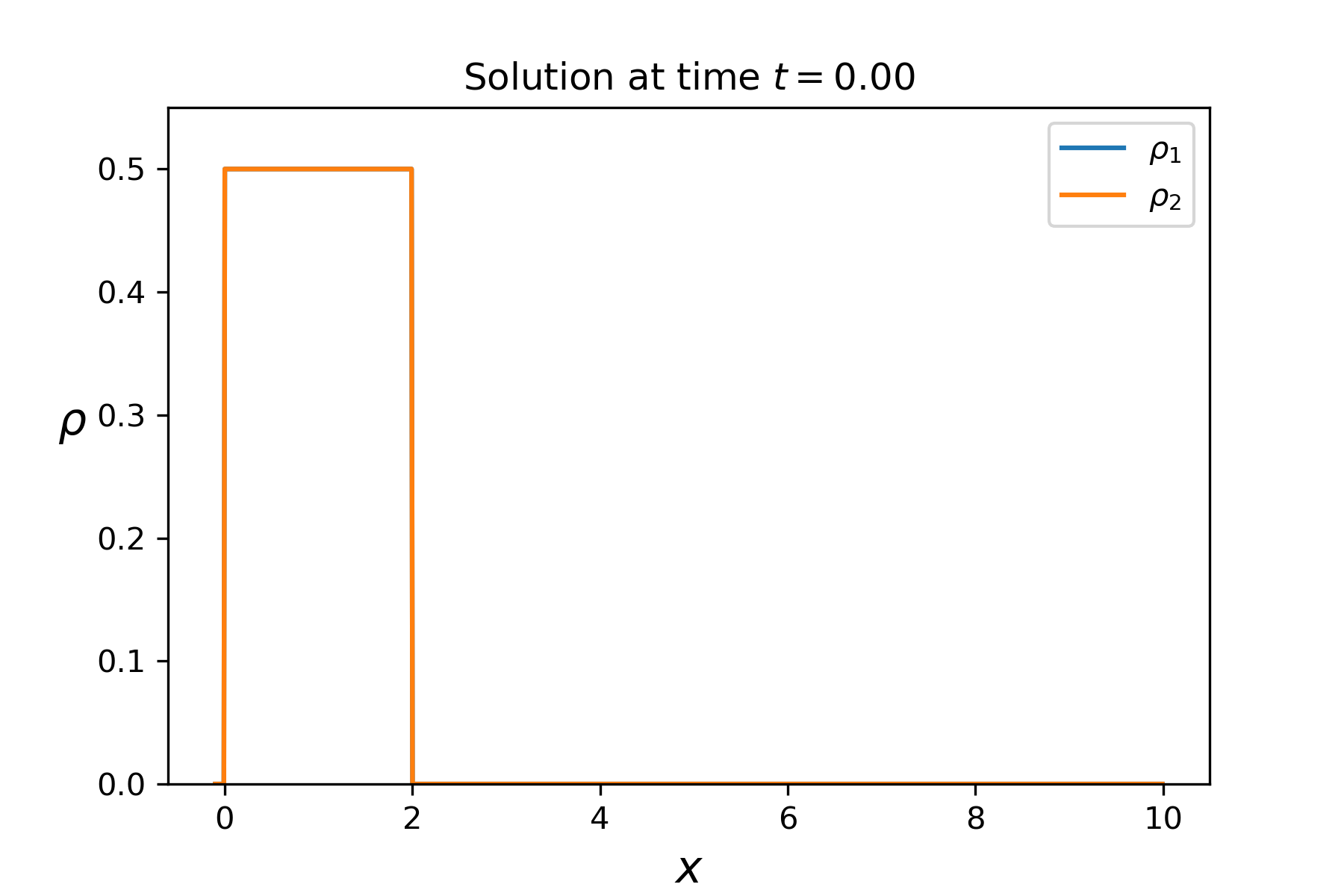}
  \end{subfigure}%
  \begin{subfigure}{0.33\textwidth}
    \centering \includegraphics[width=\textwidth, trim=20 10 20 10
    keepaspectratio]{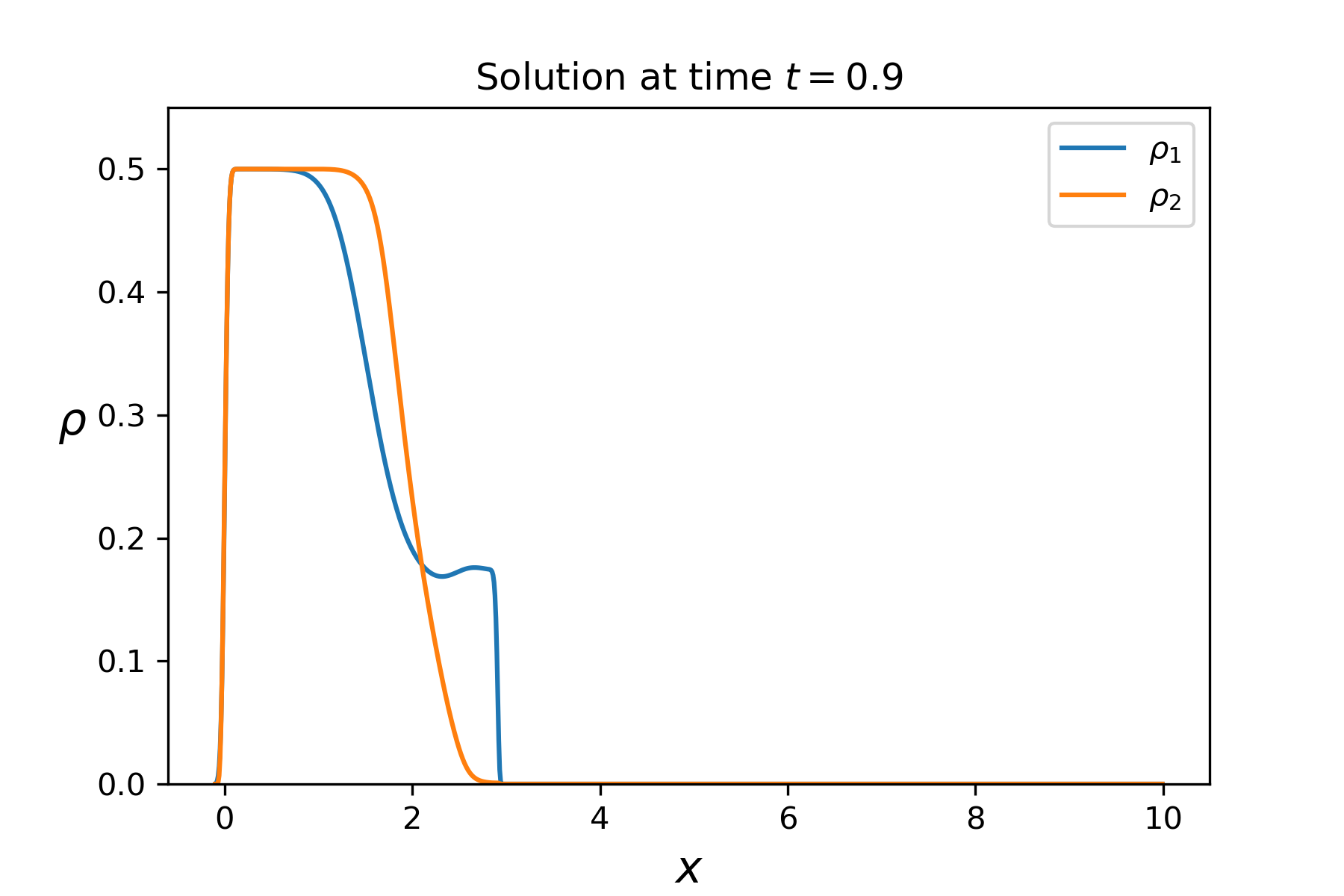}
  \end{subfigure}%
  \begin{subfigure}{0.33\textwidth}
    \centering \includegraphics[width=\textwidth, trim=20 10 20 10
    keepaspectratio]{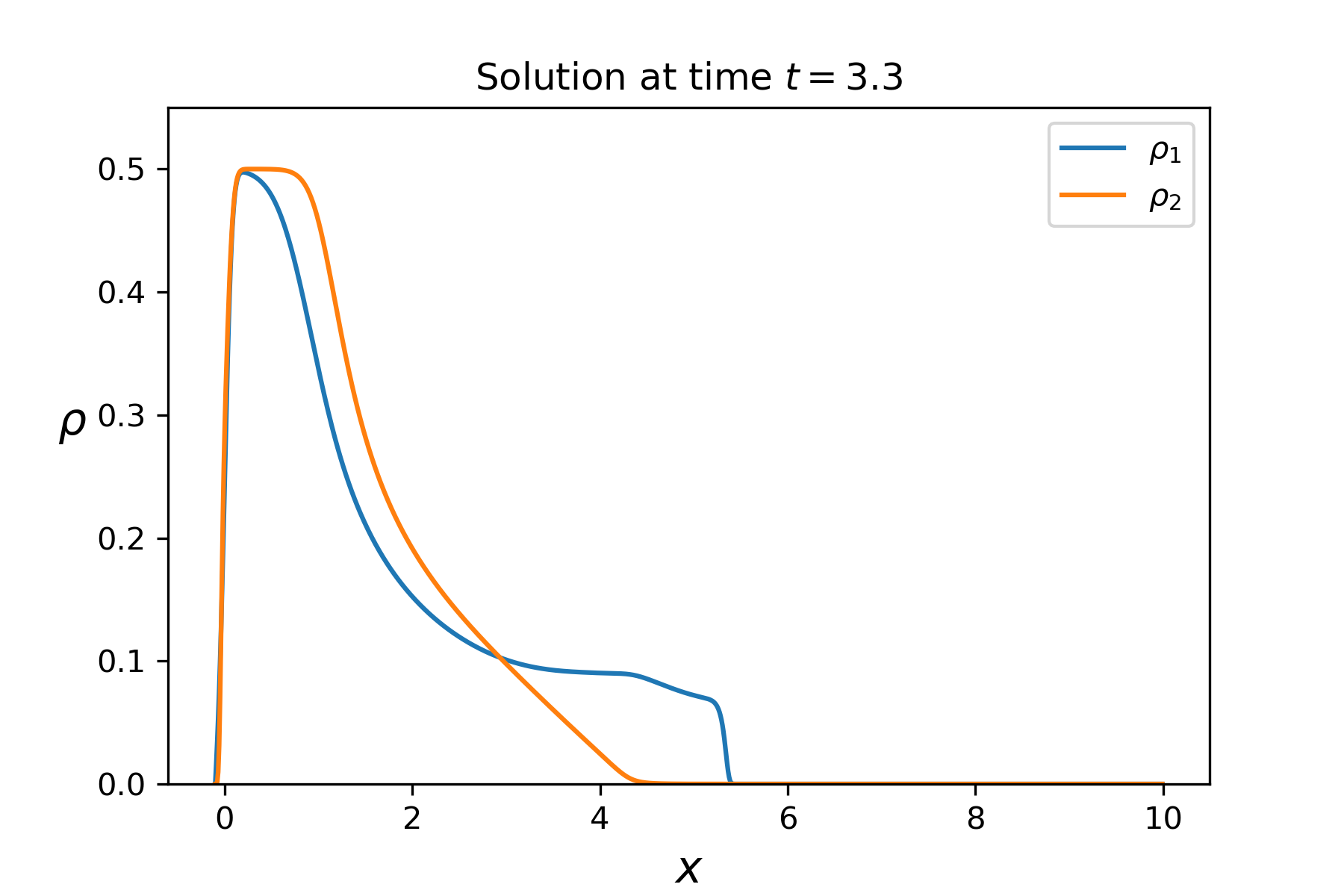}
  \end{subfigure}\\
  \begin{subfigure}{0.33\textwidth}
    \centering \includegraphics[width=\textwidth, trim=20 10 20 10
    keepaspectratio]{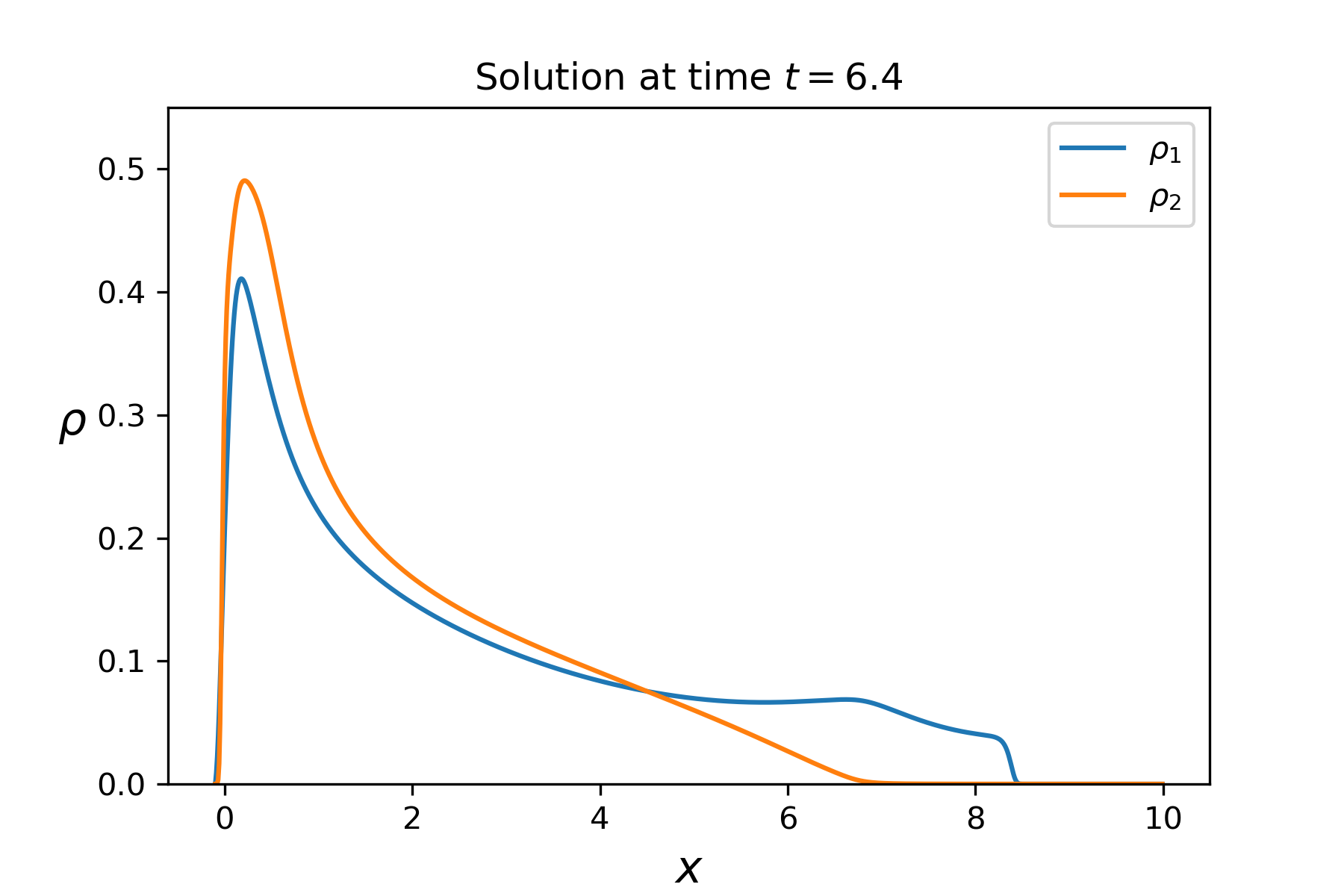}
  \end{subfigure}%
  \begin{subfigure}{0.33\textwidth}
    \centering \includegraphics[width=\textwidth, trim=20 10 20 10
    keepaspectratio]{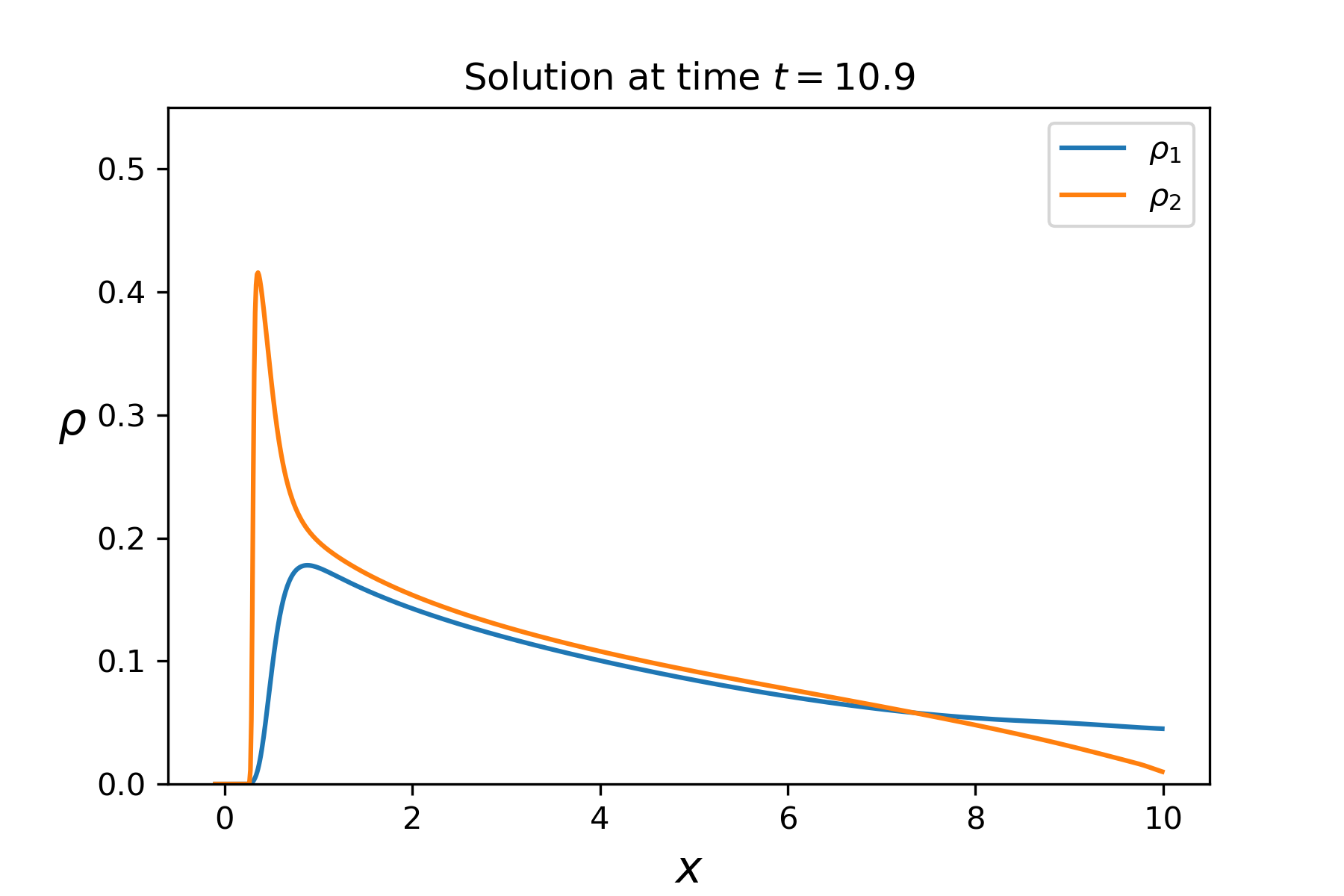}
  \end{subfigure}%
  \begin{subfigure}{0.33\textwidth}
    \centering \includegraphics[width=\textwidth, trim=20 10 20 10
    keepaspectratio]{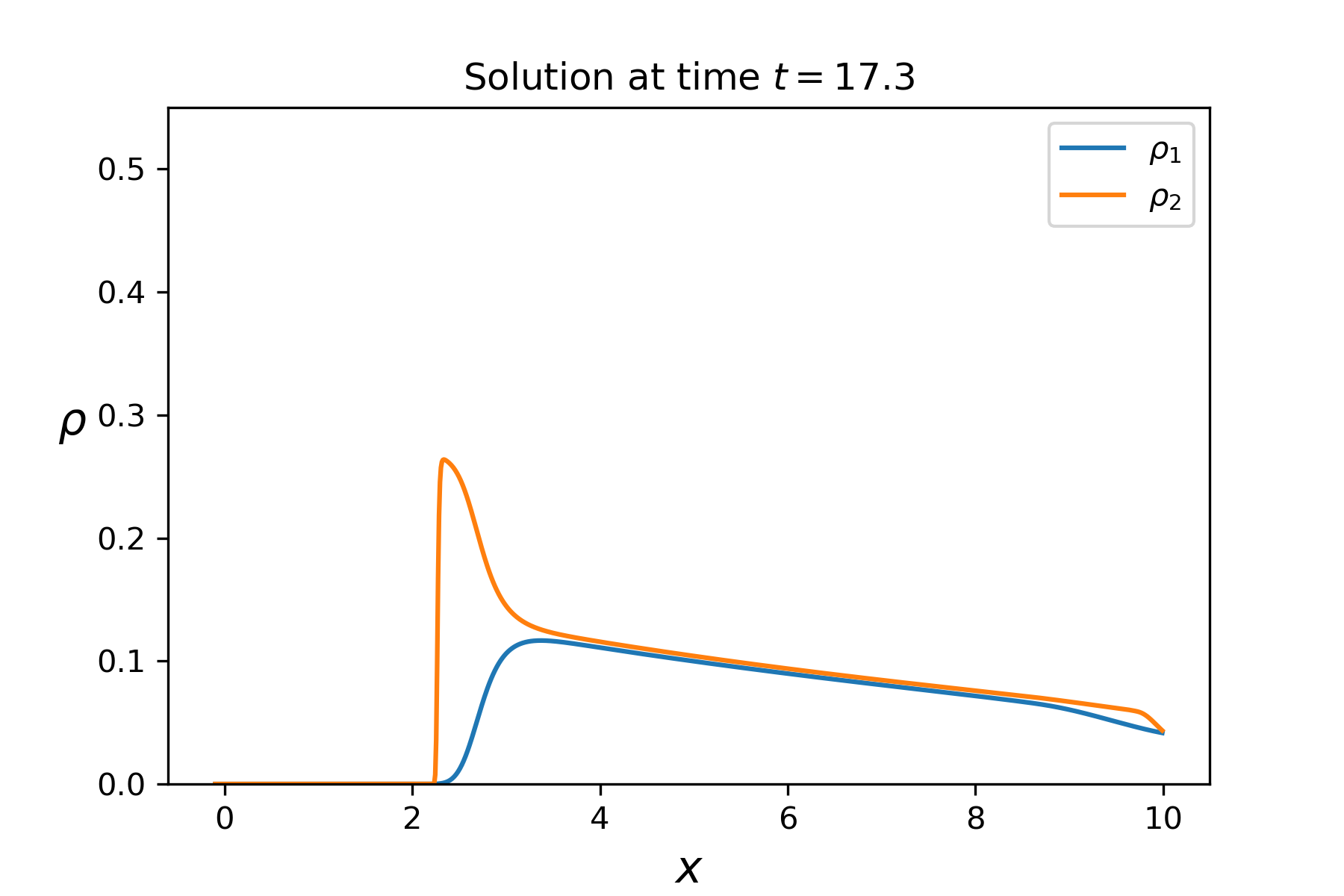}
  \end{subfigure}
  \caption{Solutions to~\eqref{eq:9} for $n=2$ with the parameters as
    in \S~\ref{subsec: horizon}. The two populations have the same
    speed law, but the forward horizon of the first population is $5$
    times that of the second one. As a result, the first population
    moves faster and its front vehicles share almost the same
    speed. On the contrary, the rightmost part of the graph of $\rho_2$
    is less steep.}
  \label{fig:horizon}
\end{figure}

We choose the following speed law, satisfying~\ref{item:1}:
\begin{equation}
  \label{eq:26}
  v_i(t, x, q)
  \coloneqq
  \left\{
  \begin{array}{lr@{\,}c@{\,}l}
    1
    & q
    & <
    & 0
    \\
    ( 1 -  q)^3
    & q
    & \in
    & [0,1]
    \\
    0
    & q
    & >
    & 1
  \end{array}
  \right.
  \qquad q = \eta_{i1}*\rho_1 + \eta_{i2}*\rho_2\,,
  \qquad i=1,2\,.
\end{equation}
System~\eqref{eq:9} is equipped with the initial datum
\begin{equation*}
  \rho^o_1(x)=\rho^o_2(x)=0.5 \, \chi_{[0,2]}(x) \,,
\end{equation*}
as in Figure~\ref{fig:horizon}, top left.
The resulting evolution, approximated by a numerical integration on a
mesh of $10000$ points on the numerical domain $[0,\, 10]$, shows the
relevance of the forward horizon, see Figure~\ref{fig:horizon}. For
the first population, the free space on the road stretch ahead has
more relevance in the choice of the speed. As a consequence, the front
vehicles in the first population are faster. Moreover, again due to
the length of the forward horizon, the front vehicles almost share the
same speed, as shown in the graphs at times $t=0.9, \, 3.3, \, 6.4$,
where the density of the first population displays a rather steep
right front. On the contrary, the slope in the right part of the graph
of $\rho_2$ is significantly lower.

\subsection{Different Maximal Speeds}
\label{subsec:overtaking}

On the numerical domain $[0,\,100]$ we now consider $n=3$ populations
differing only in their maximal speed, i.e., we set for $i=1,2,3$
\begin{equation}
  \label{eq:30}
  v_i(t, x, q)
  \coloneqq
  \left\{
    \begin{array}{lr@{\,}c@{\,}l}
      V_i
      & q
      & <
      & 0
      \\
      V_i( 1 - q )^3
      & q
      & \in
      & [0,1]
      \\
      0
      & q
      & >
      & 1
    \end{array}
  \right.
  \quad
  q = \eta_{i1}* \rho_1 + \eta_{i2}* \rho_2+ \eta_{i3} * \rho_3\,,
  \quad
  \begin{array}{rcl}
    V_1
    & =
    & 1.5
    \\
    V_2
    & =
    & 0.9
    \\
    V_3
    & =
    & 0.5 \,.
  \end{array}
\end{equation}
We choose the initial datum
\begin{equation}
  \label{eq:29}
  (\rho_o)_1(x) = 0.3 \, \chi_{[1,5]}(x)
  \,,\quad
  (\rho_o)_2(x) = 0.3 \, \chi_{[8,12]}(x)
  \,,\quad
  (\rho_o)_3(x) = 0.3 \, \chi_{[15,19]}(x)
\end{equation}
and kernel functions~\eqref{eq:kernel} with forward horizon $f_{ij}=1.0$ and backward horizon $b_{ij}=0.01$ for $i,j=1,2,3$.

\begin{figure}[!ht]
  \begin{subfigure}{0.33\textwidth}
    \centering \includegraphics[width=1\textwidth, trim=20 0 20 0
    keepaspectratio]{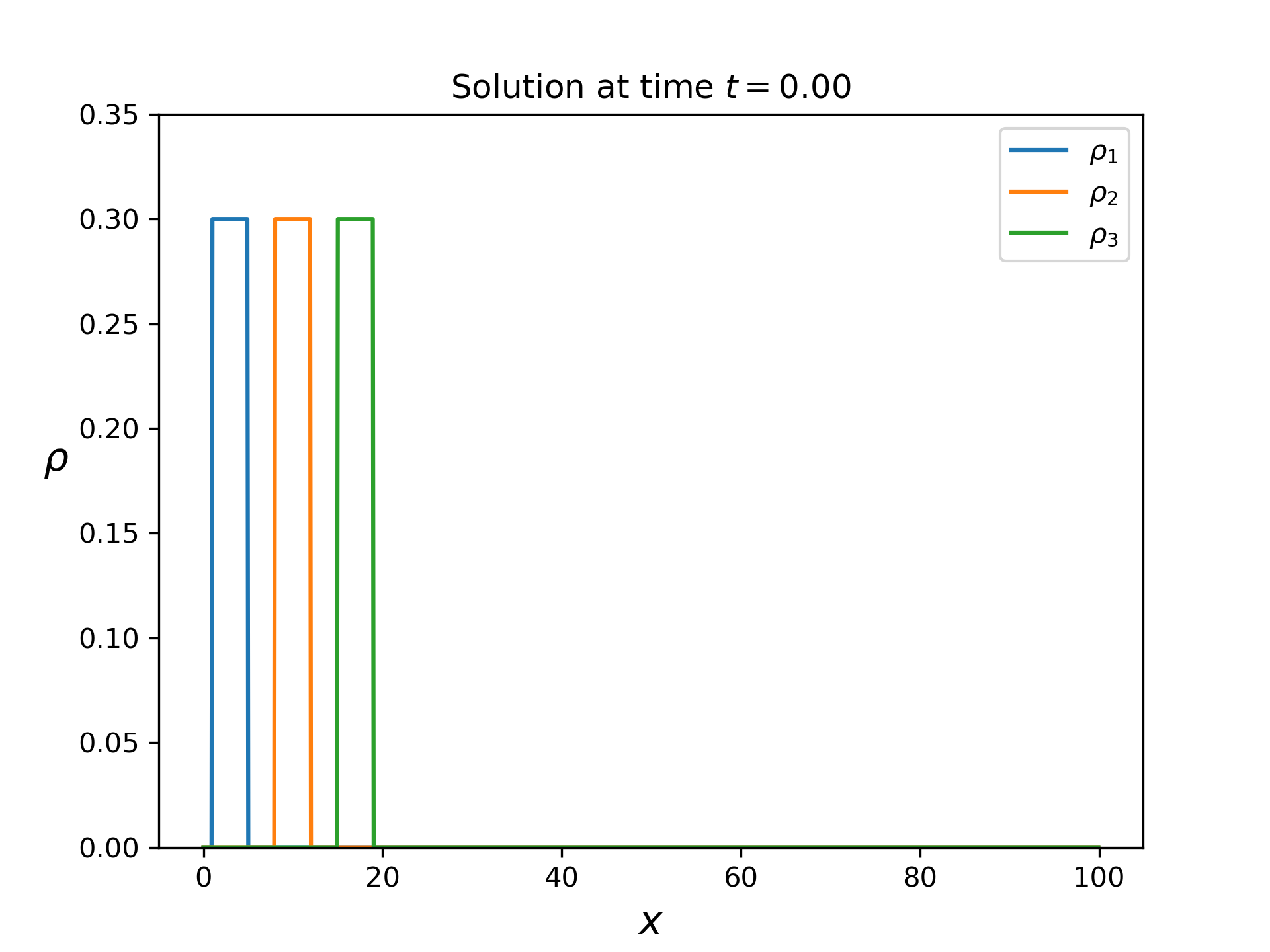}
  \end{subfigure}%
  \begin{subfigure}{0.33\textwidth}
    \centering \includegraphics[width=\textwidth, trim=20 0 20 0
    keepaspectratio]{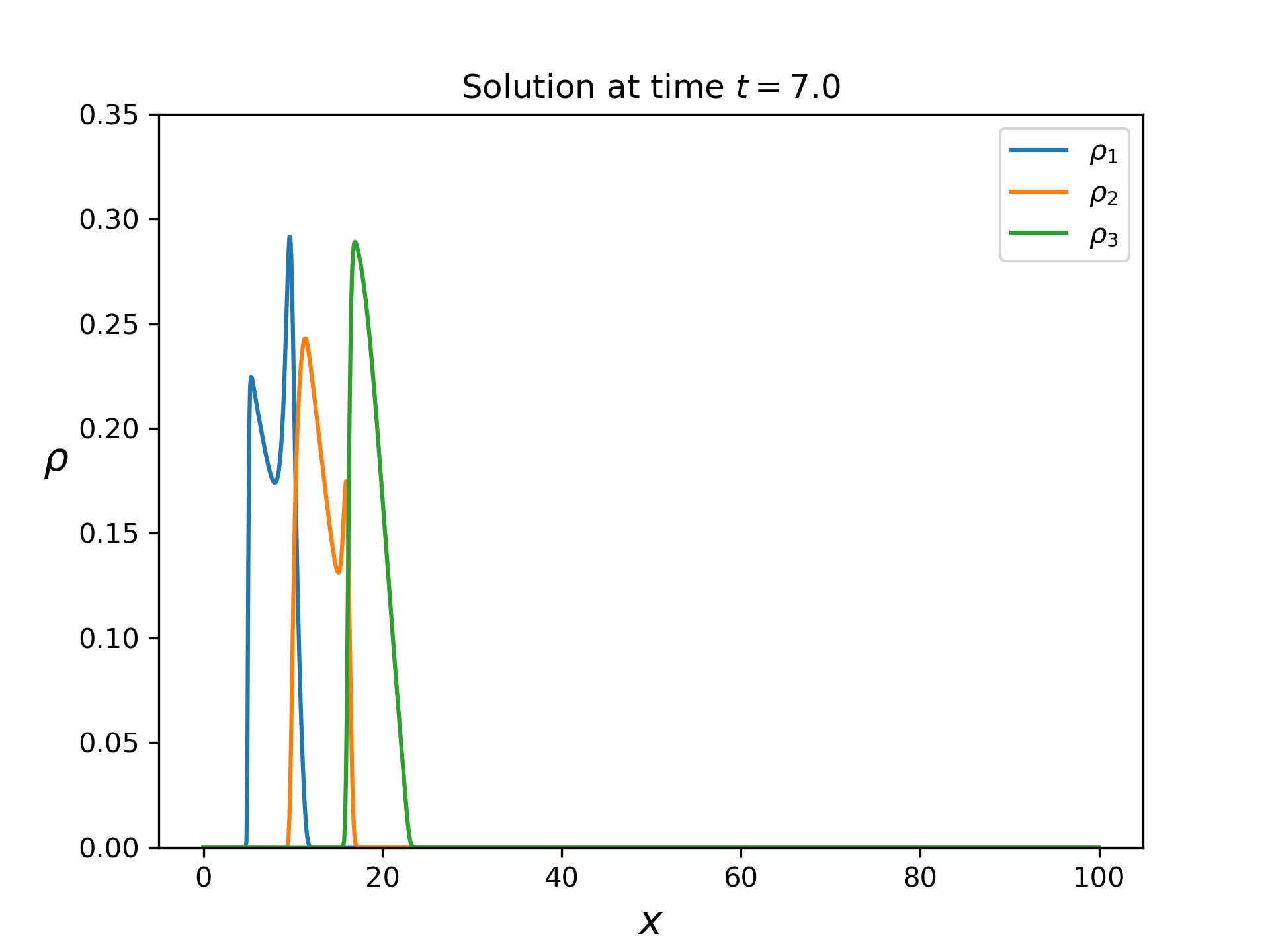}
  \end{subfigure}%
  \begin{subfigure}{0.33\textwidth}
    \centering \includegraphics[width=\textwidth, trim=20 0 20 0
    keepaspectratio]{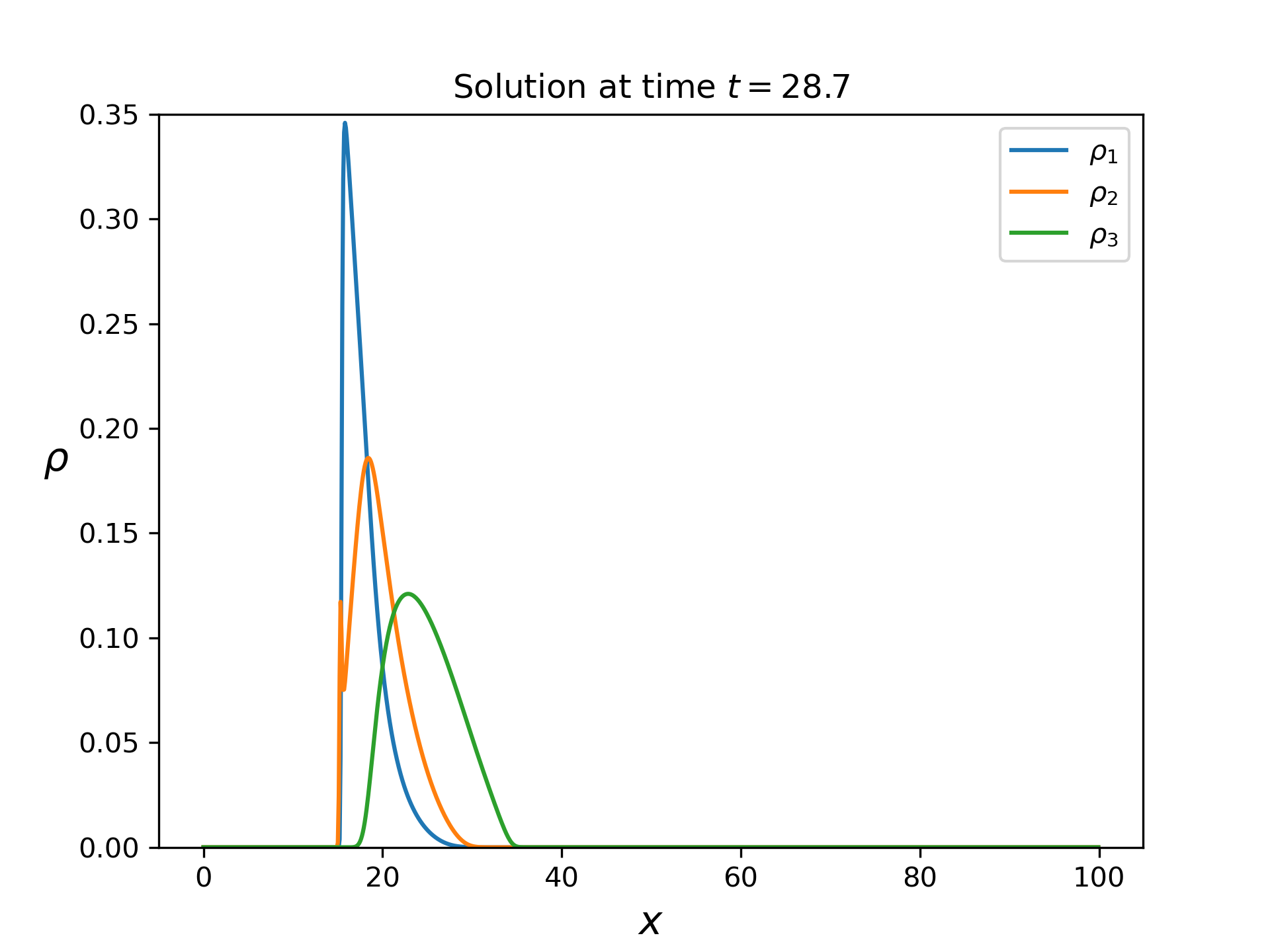}
  \end{subfigure}\\
  \begin{subfigure}{0.33\textwidth}
    \centering \includegraphics[width=\textwidth, trim=20 0 20 0
    keepaspectratio]{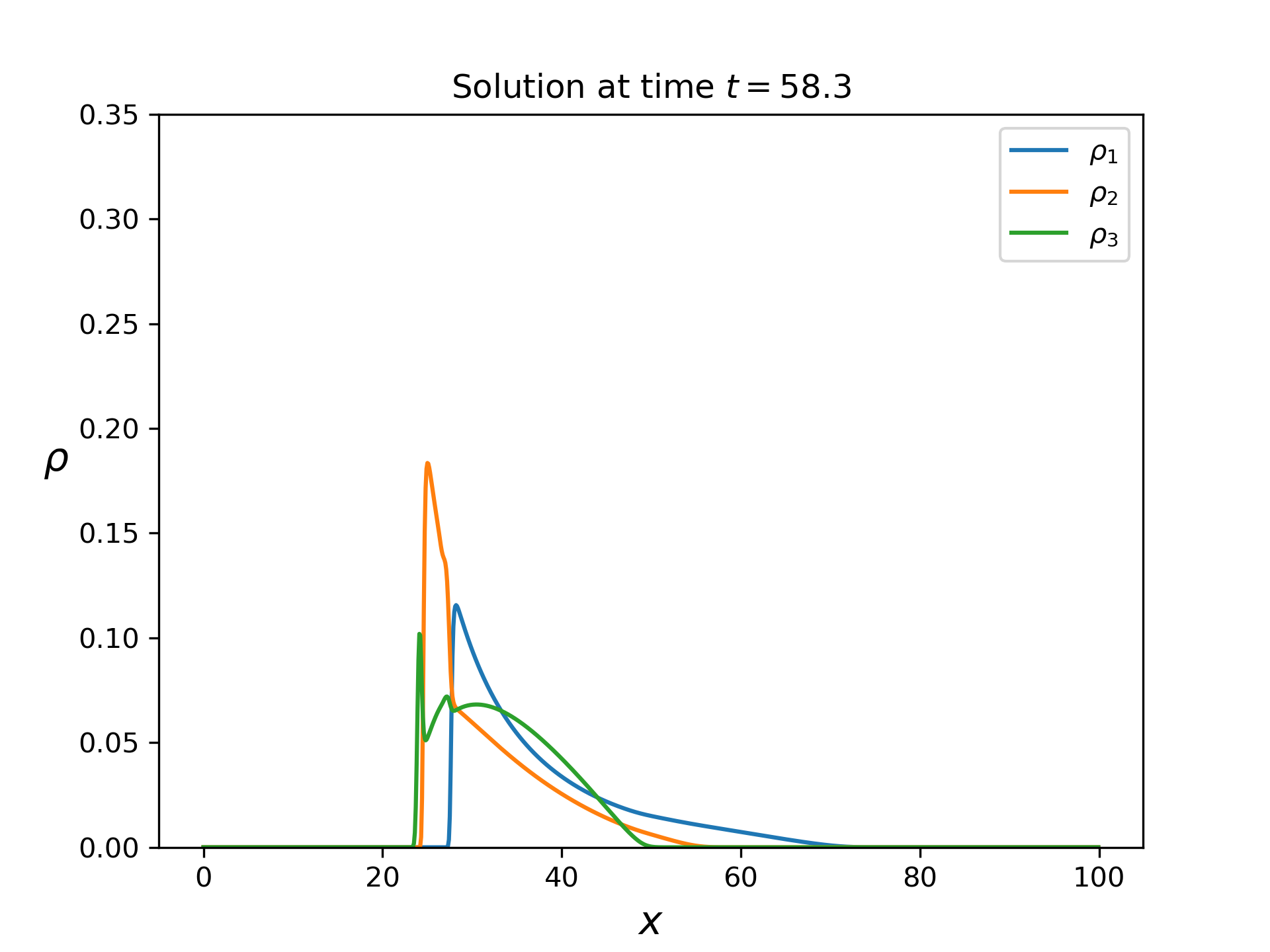}
  \end{subfigure}%
  \begin{subfigure}{0.33\textwidth}
    \centering \includegraphics[width=\textwidth, trim=20 0 20 0
    keepaspectratio]{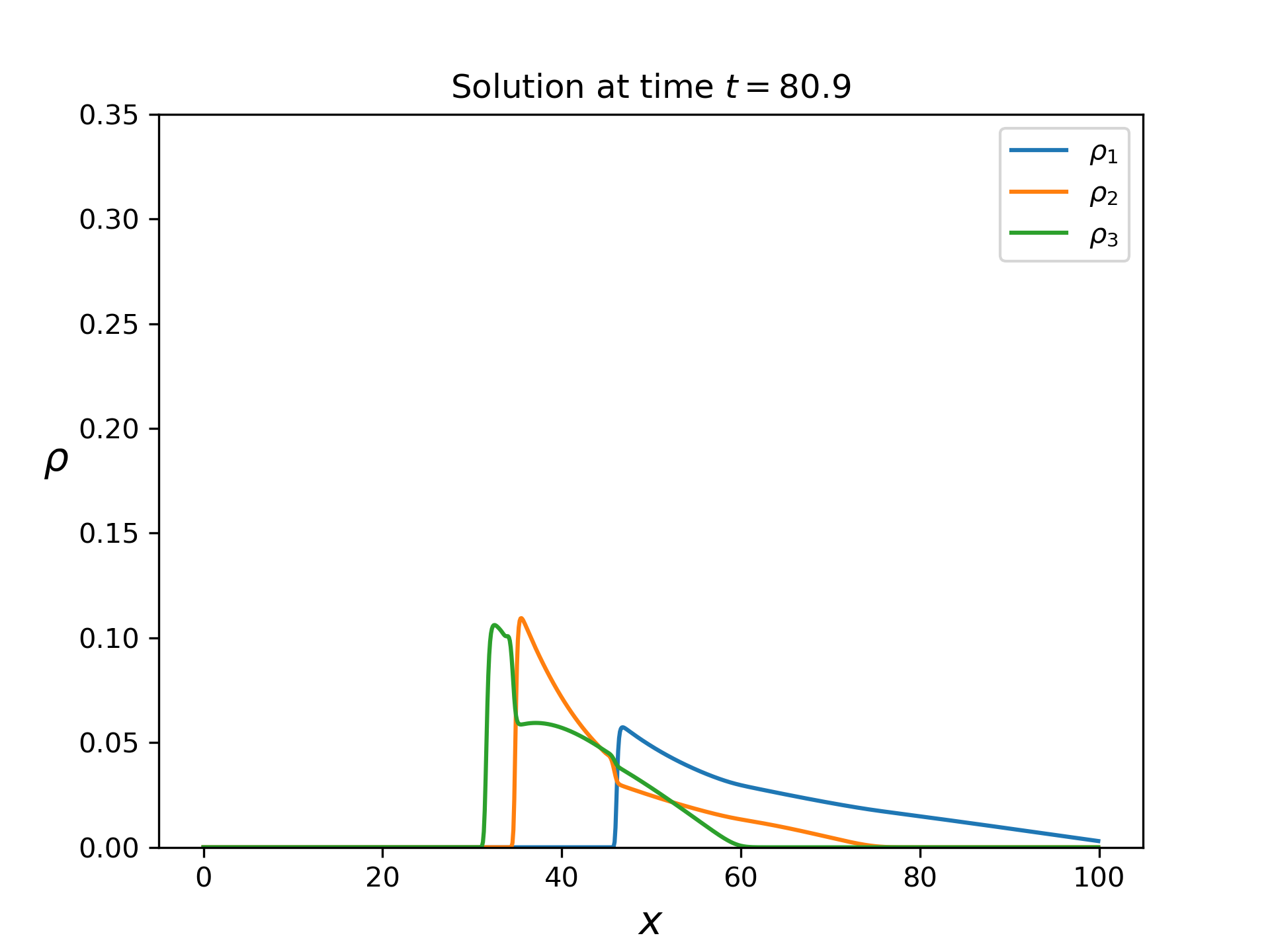}
  \end{subfigure}%
  \begin{subfigure}{0.33\textwidth}
    \centering \includegraphics[width=\textwidth, trim=20 0 20 0
    keepaspectratio]{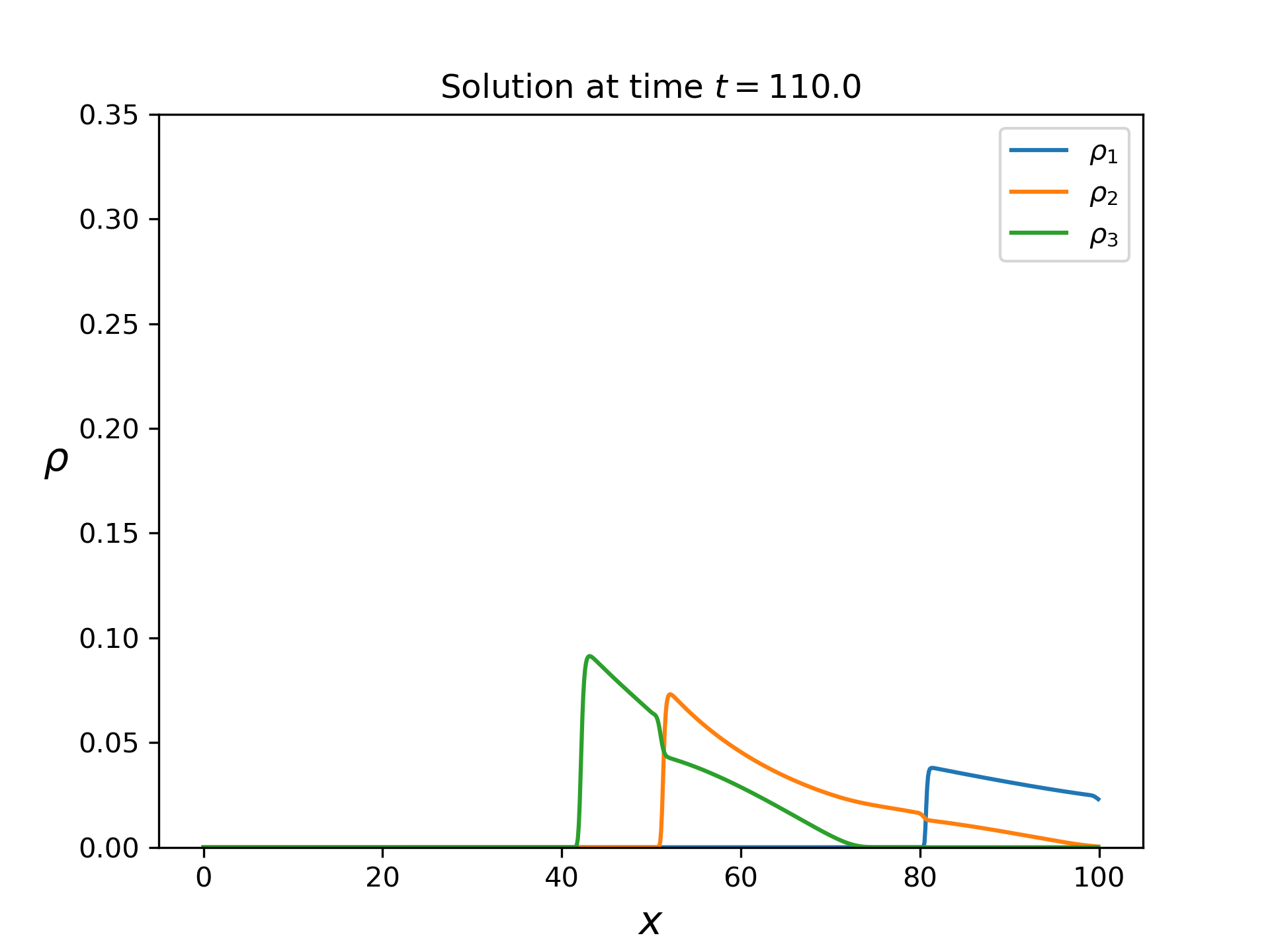}
  \end{subfigure}
  \caption{Solutions to~\eqref{eq:9}--\eqref{eq:30}--\eqref{eq:29}
    with $n=3$ populations on the domain $[0,\, 100]$ with the
    parameters given in Subsection~\ref{subsec:overtaking}. Note that
    faster populations overtake slower populations and towards the end
    of the integration the populations appear ordered according to
    their maximal speed. Moreover, the graph at $t=28.7$ shows that a
    density may reach values higher than those attained at the initial
    time.}
  \label{fig:sorpasso}
\end{figure}

Figure~\ref{fig:sorpasso} shows the result of the numerical
integration of~\eqref{eq:9}--\eqref{eq:30}--\eqref{eq:29} obtained
with a uniform mesh of $10000$ points.
At the initial time, the slowest population is ahead and the fastest
is behind. Model~\eqref{eq:9} allows overtakes and, in fact, at the
final time the populations are ordered according to their maximal
speeds. It is realistic that during overtakes queues form, see in
particular the graphs at times $t=7.0,\, 28.7,\,80.9$ in
Figure~\ref{fig:sorpasso}. Note that the graph at time $t=28.7$
clearly shows that the density $\rho_1 (t)$ exceeds the maximum of the
initial datum $(\rho_o)_1$, coherently with the lack of validity of
the maximum principle in the present setting.

\subsection{Space Inhomogeneity and Forward Horizon}
\label{subsec:space-inhom-forw}

The generality of system~\eqref{eq:9} allows to evaluate the combined
effects of space inhomogeneities and nonlocal terms. The numerical
domain $[0,\, 20]$ is a road stretch with a bottleneck in $[5,\, 10]$
where the maximal speed smoothly decreases up to $50\%$. We show below
that vehicles with a larger forward horizon better cope with the
slowdown and pass through it faster than vehicles with a lower
one. More precisely, we confront a solution to~\eqref{eq:9} with that
to the standard Lighthill-Whitham~\cite{LighthillWhitham} and
Richards~\cite{Richards} model.

First, consider~\eqref{eq:9} with $n=1$, $\eta$ as in~\eqref{eq:kernel} with $f_{11} = 1.0$, $b_{11} = 0.01$ and
\begin{eqnarray}
  \label{eq:22}
  v(t, x, q)
  & \coloneqq
  & \left\{
    \begin{array}{lr@{\,}c@{\,}l}
      \rev{V}(x)
      & q
      & <
      &0
      \\
      \rev{V}(x) \; ( 1 -  q)^3
      & q
      & \in
      & [0,1]
      \\
      0
      & q
      & >
      & 1
    \end{array}
        \right.
        \qquad q = \eta * \rho\,,
  \\
  \label{eq:23}
  \rev{V}(x)
  & =
  & \begin{cases}
    1 - \frac{32}{5^6} \, (x-5)^3 \, (10-x)^3
    & \text{if}\quad x \in [5, 10]
    \\
    1
    &\quad \text{otherwise},
  \end{cases}
  \\
  \label{eq:27}
  \rho_o (x)
  & =
  & 0.8 \, \caratt{[1,\,3]} (x)
\end{eqnarray}
so that Theorem~\ref{thm:main} applies. We compare the resulting
solution $\rho$ with that, say $r$, to the LWR model
\begin{equation}
  \label{eq:28}
  \left\{
    \begin{array}{l}
      \partial_t r + \partial_x \left(r \, v (\rev{t}, x,r)\right) = 0
      \\
      r (0,x) = \rho_o (x)
    \end{array}
  \right.
\end{equation}
where $v$ is as in~\eqref{eq:22}--\eqref{eq:23} and $\rho_o$ is as
in~\eqref{eq:27}. Both solutions are computed on a mesh of $10000$
points.
\begin{figure}[!ht]
  \begin{subfigure}{0.33\textwidth}
    \centering \includegraphics[width=1\textwidth, trim=20 0 20 0
    keepaspectratio]{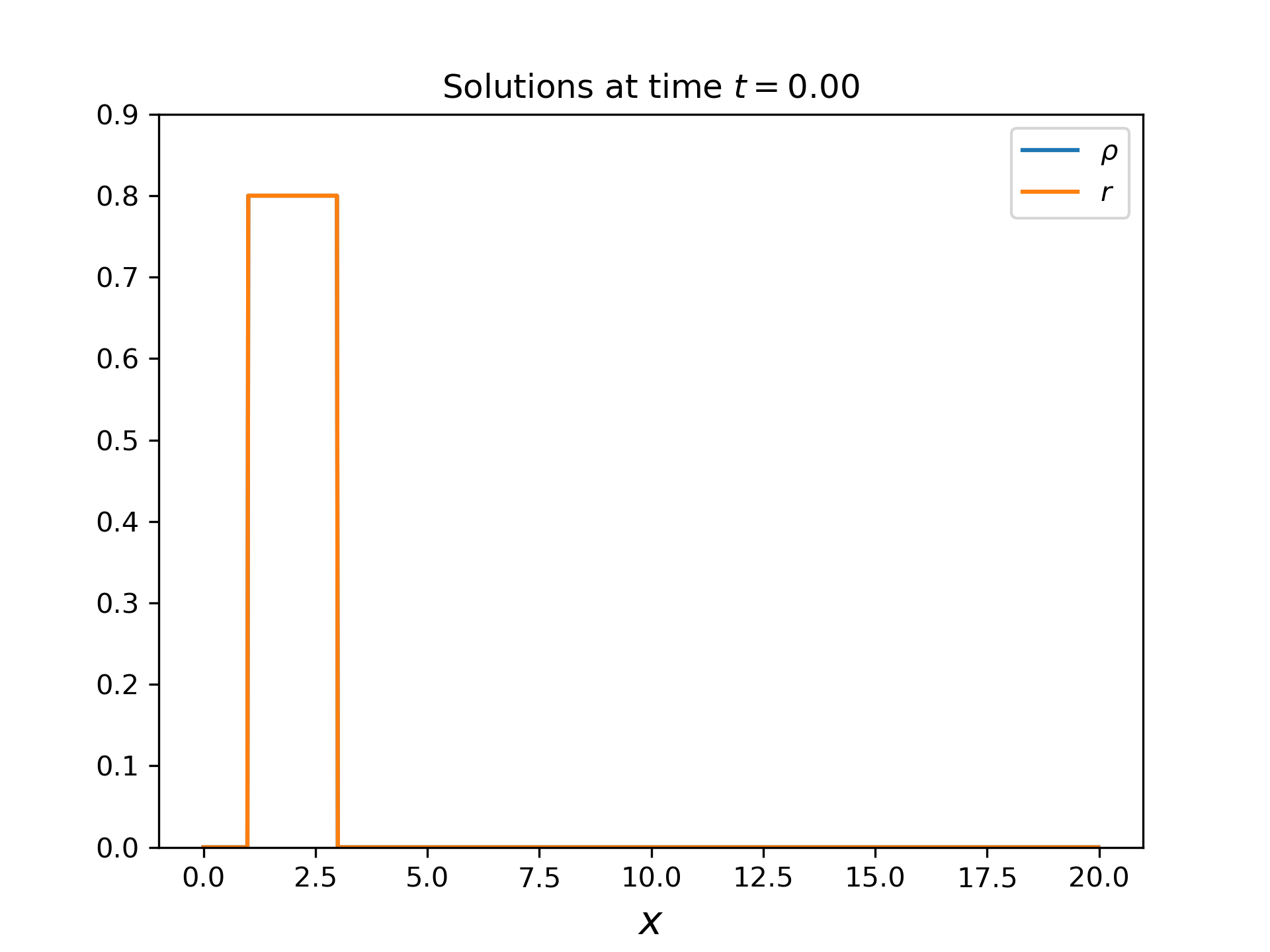}
  \end{subfigure}%
  \begin{subfigure}{0.33\textwidth}
    \centering \includegraphics[width=\textwidth, trim=20 0 20 0
    keepaspectratio]{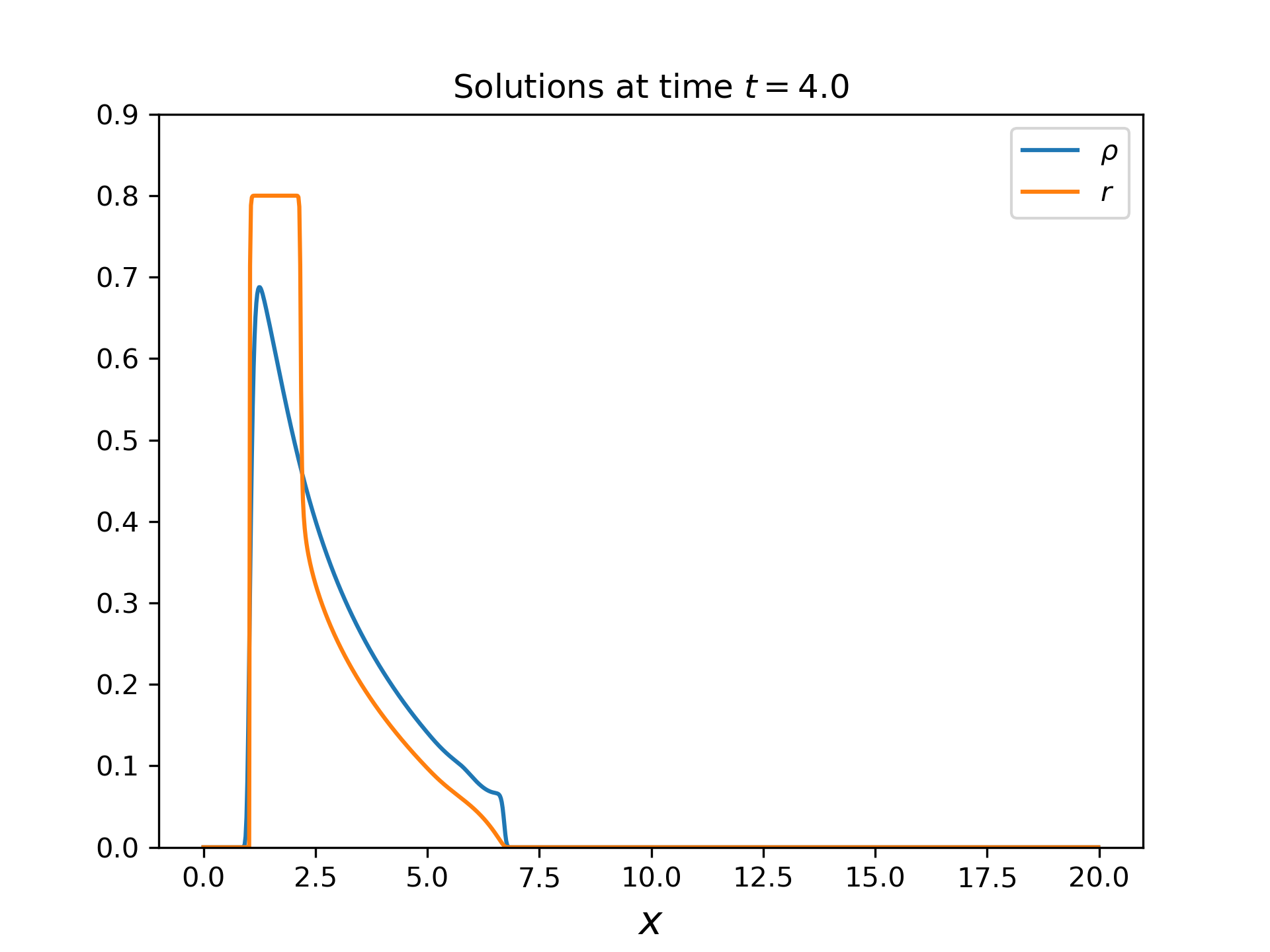}
  \end{subfigure}%
  \begin{subfigure}{0.33\textwidth}
    \centering \includegraphics[width=\textwidth, trim=20 0 20 0
    keepaspectratio]{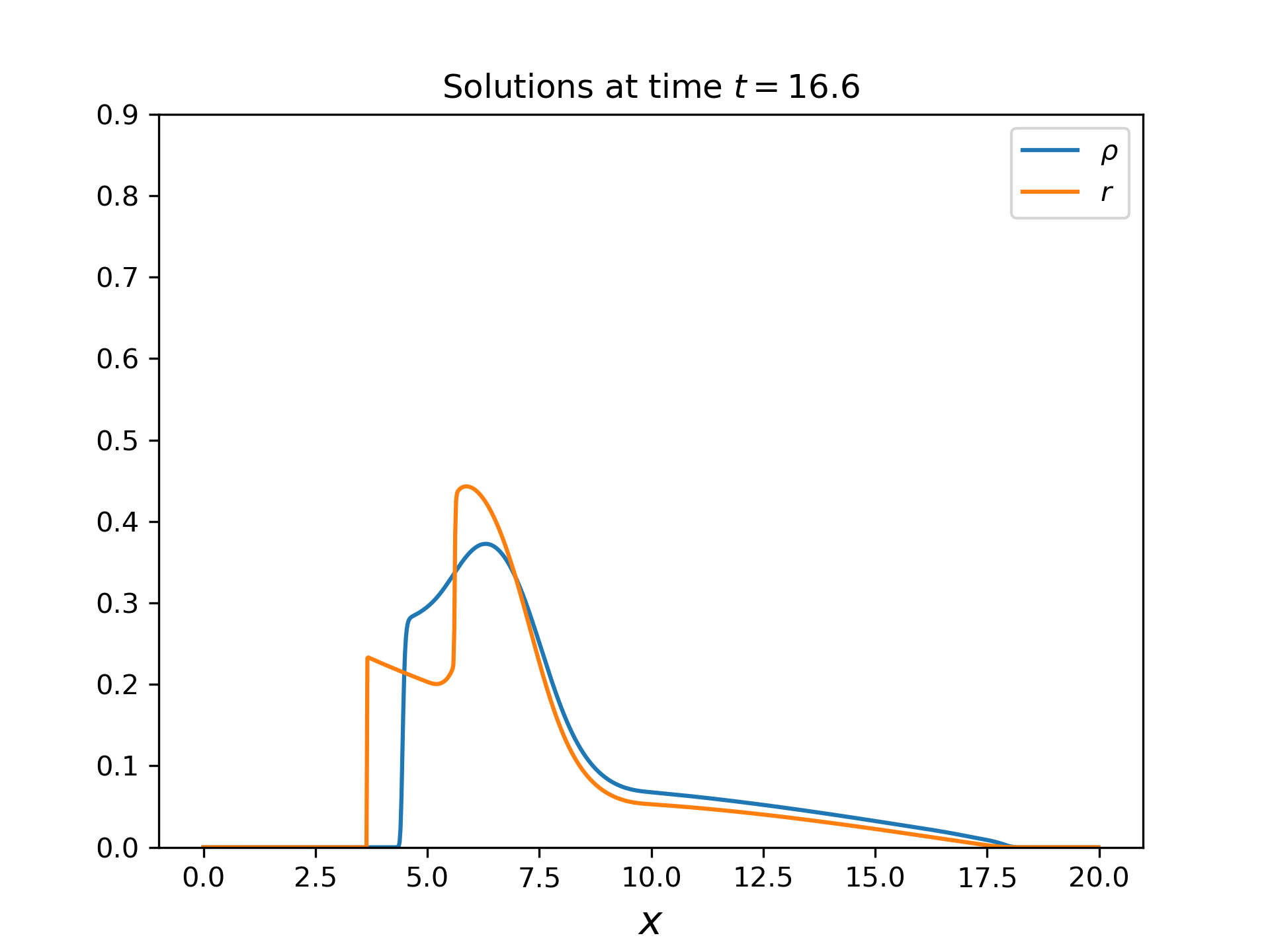}
  \end{subfigure}\\
  \begin{subfigure}{0.33\textwidth}
    \centering \includegraphics[width=\textwidth, trim=20 0 20 0
    keepaspectratio]{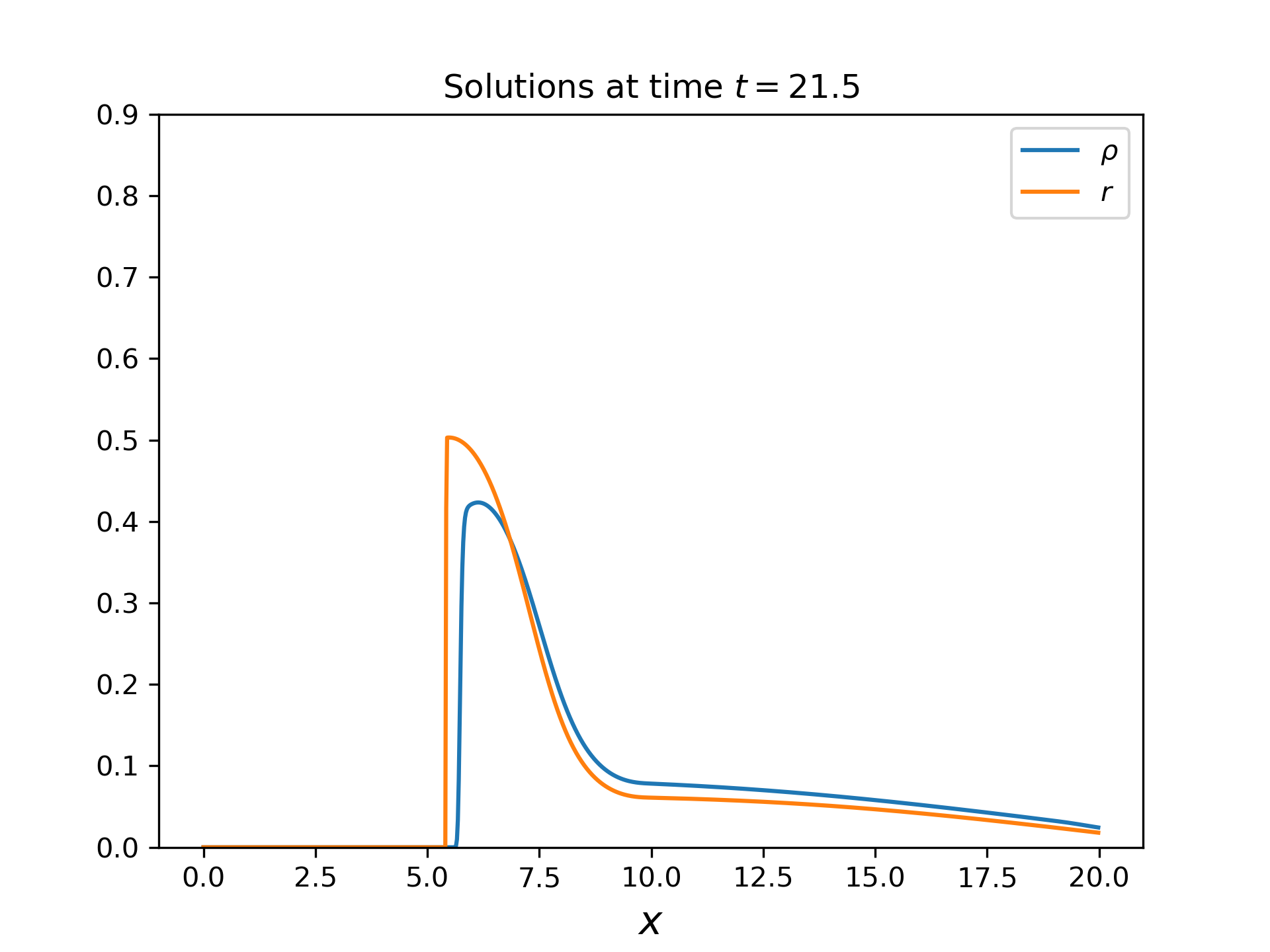}
  \end{subfigure}%
  \begin{subfigure}{0.33\textwidth}
    \centering \includegraphics[width=\textwidth, trim=20 0 20 0
    keepaspectratio]{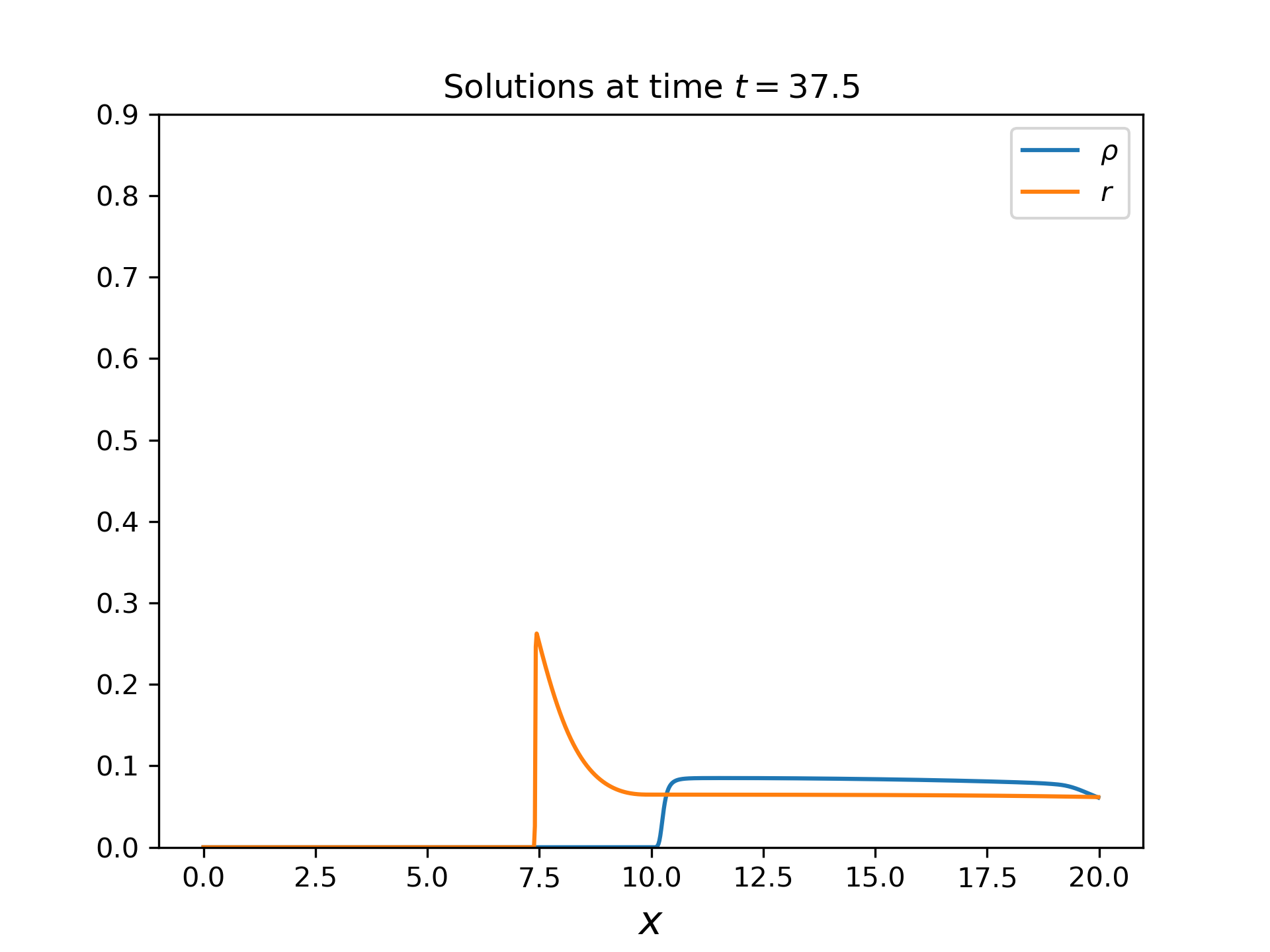}
  \end{subfigure}%
  \begin{subfigure}{0.33\textwidth}
    \centering \includegraphics[width=\textwidth, trim=20 0 20 0
    keepaspectratio]{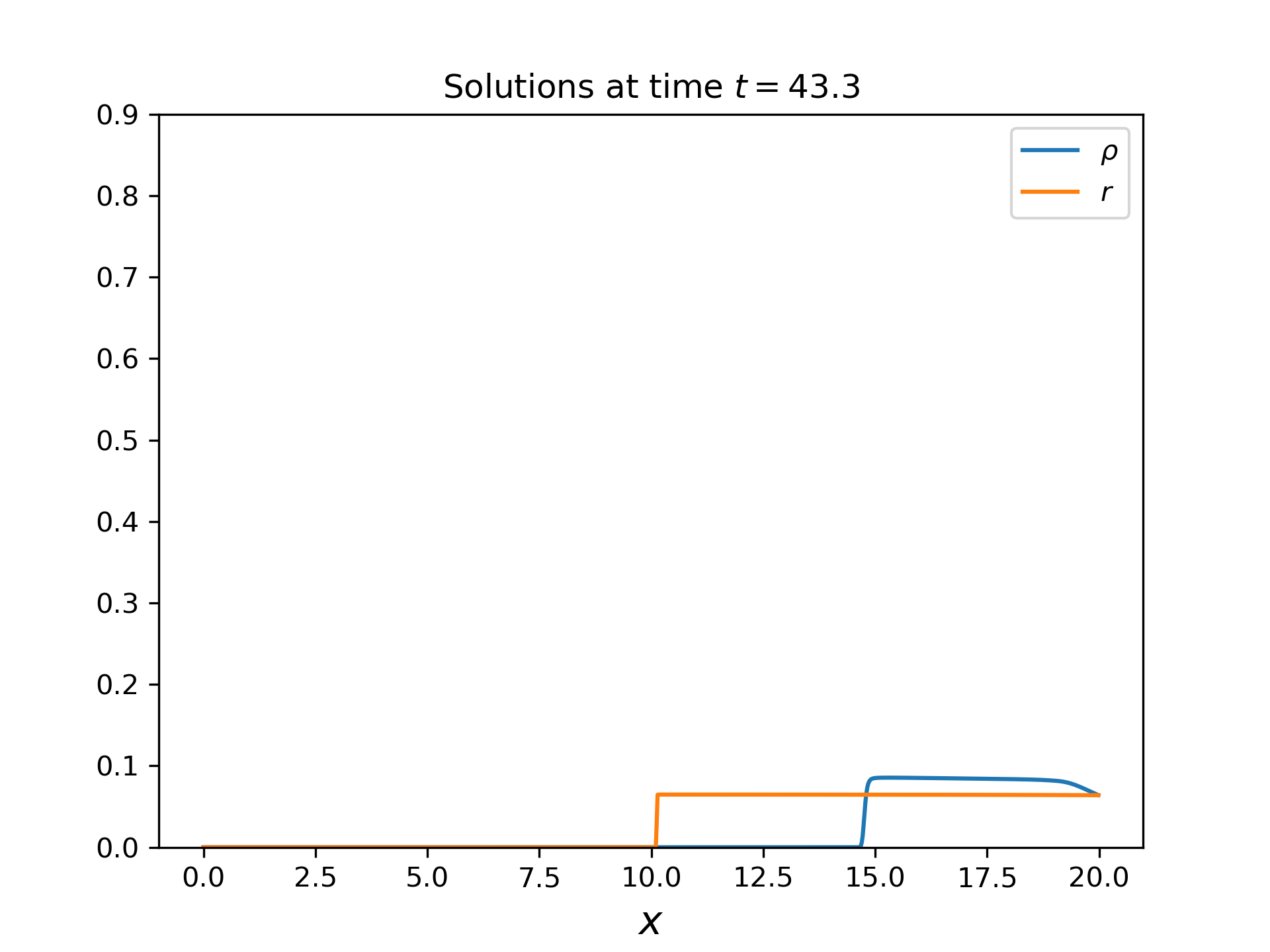}
  \end{subfigure}
  \caption{Solution $\rho$ to the nonlocal
    equation~\eqref{eq:9}--\eqref{eq:22}--\eqref{eq:23}--\eqref{eq:27}
    compared to $r$, solving the LWR model~\eqref{eq:28}, with the
    same initial datum~\eqref{eq:27} and speed
    law~\eqref{eq:22}--\eqref{eq:23}. The effect of the bottleneck in
    $[5,\,10]$ is clearly seen since time $t=4.0$. Note that already
    at time $t=37.5$ the solution to the nonlocal problem is supported
    on the right of the bottleneck, while the same happens only after
    $t=43.3$ for the solution to the LWR model.}
  \label{fig:street}
\end{figure}

Figure~\ref{fig:street} displays the solutions $\rho$ and $r$ to the
two problems. It stems out that the vehicles with a positive forward
horizon slow down well before the bottleneck, cause less congestion
and pass through the bottleneck in less time than the others. Indeed,
at time $t=37.5$ the solution to~\eqref{eq:9} is supported on the
right of the bottleneck, while only after $t=43.3$ the solution
to~\eqref{eq:28} has overcome the slowdown.

\section{Analytical Proofs}
\label{sec:technical-details}

Without loss of generality, we assume that $t_o = 0$.
Throughout, we use the Euclidean norm $\norma{\,\cdot\,}$ in
$\reali^n$. As usual, we equip the spaces $\L1 (\reali; \reali^n)$,
$\W1\infty (\reali; \reali^n)$, and $\W2\infty (\reali; \reali^n)$
with the norms
\begin{displaymath}
  \begin{array}{@{}rcl@{\qquad}rcl@{}}
    \norma{\rho}_{\L1 (\reali;\reali^n)}
    & \coloneqq
    &\int_{\reali} \norma{\rho (x)} \d{x}
    &
      \norma{\rho}_{\W1\infty (\reali;\reali^n)}
    & \coloneqq
    & \norma{\rho}_{\L\infty (\reali;\reali^n)} + \norma{\frac{\d\rho}{\d{x}}}_{\L\infty (\reali;\reali^n)}
    \\
    \norma{\rho}_{\L\infty (\reali;\reali^n)}
    & \coloneqq
    &\esssup_{x \in \reali} \norma{\rho (x)}
    &
      \norma{\rho}_{\W2\infty (\reali;\reali^n)}
    & \coloneqq
    & \norma{\rho}_{\W1\infty (\reali;\reali^n)} + \norma{\frac{\mathinner{\mathrm{d}^2{\rho}}}{\d{x^2}}}_{\L\infty (\reali;\reali^n)} \,.
  \end{array}
\end{displaymath}
By $\tv (u)$ we denote the total variation of $u$ on its domain, while
$\tv (u,I)$ is the total variation of $u$ on the real interval $I$,
see~\cite[Chapter~3]{MR1857292}.

The following basic property of $\BV$ functions is a slight extension
of~\cite[Lemma~2.3]{MR1816648}.

\begin{lemma}
  \label{lem:TV}
  For any $u \in \BV (\reali; \reali)$ and
  $h \in \L\infty (\reali; \reali)$,
  \begin{displaymath}
    \int_{\reali} \modulo{u\left(x+h (x)\right) - u (x)} \d{x}
    \leq
    \left( \max\left\{ \esssup\nolimits_{\reali} h, 0\right\}
      - \min\left\{ \essinf\nolimits_{\reali} h, 0\right\}\right)
    \tv (u) \,.
  \end{displaymath}
\end{lemma}

\begin{proof}
  Let $U(x) \coloneqq \tv\left(u, \mathopen]-\infty,
    x\mathclose[\right)$. Then, by~\cite[Lemma~2.3]{MR1816648}
  \begin{eqnarray*}
    &
    & \int_{\reali} \modulo{u\left(x+h(x)\right) -u(x)} \d{x}
    \\
    & \leq
    & \int_{\reali}
      \left(
      U\left(x+ \max\left\{ \esssup\nolimits_{\reali} h, 0\right\} \right)  -
      U\left( x + \min\left\{ \essinf\nolimits_{\reali} h, 0\right\}\right) \right)
      \,\d{x}
    \\
    & =
    & \int_{\reali}
      \left(U\left(y + \left(\max\left\{ \esssup\nolimits_{\reali} h, 0\right\} - \min\left\{ \essinf\nolimits_{\reali} h, 0\right\}\right)  \right) - U(y)\right)
      \d{y}
    \\
    & \leq
    & \left( \max\left\{ \esssup\nolimits_{\reali} h, 0\right\}
      - \min\left\{ \essinf\nolimits_{\reali} h, 0\right\}\right) \tv(U)
    \\
    & =
    & \left( \max\left\{ \esssup\nolimits_{\reali} h, 0\right\}
      - \min\left\{ \essinf\nolimits_{\reali} h, 0\right\}\right) \tv(u) \,,
  \end{eqnarray*}
  completing the proof.
\end{proof}

We now recall a basic stability result on ODE flows.

\begin{proposition}
  \label{prop 2}
  Fix $T>0$ and
  $w,z \in \C0\left([0,T];\W1\infty (\reali; \reali)\right)$. For
  $t,t_o \in [0,T]$ and $x_o \in \reali$, let $t \mapsto X(t;t_o,x_o)$
  and $t \mapsto Y(t;t_o,x_o)$ be the solutions to
  \begin{equation}
    \begin{cases}
      \dot{x} = w (t, x)
      \\
      x(t_o)=x_o
    \end{cases}
    \quad \mbox{ and } \quad
    \begin{cases}
      \dot{y} = z(t,y)
      \\
      y(t_o)=x_o.
    \end{cases}
  \end{equation}
  Then for all $t \in [0,T]$,
  \begin{eqnarray*}
    \modulo{X(t;t_o,x_o)-Y(t;t_o,x_o)}
    & \leq
    & \norma{w-z}_{\C0([0,T];\L\infty(\reali;\reali))} \; \modulo{t-t_o}
    \\
    &
    & \times
      \exp \left(
      \norma{\partial_x w}_{\C0([0,T];\L\infty(\reali;\reali))} \,
      \modulo{t-t_o}
      \right) \,.
  \end{eqnarray*}
\end{proposition}

\noindent The proof is a standard consequence of Gronwall Lemma and
basic ODE theory, see~\cite[Chapter~3]{MR2347697}.

\begin{lemma}
  \label{lem:Pitilde}
  Let~\ref{item:2} hold and $v \in \mathcal{V}^n$. Call
  \begin{equation}
    \label{eq:2}
    \tilde X_M
    \coloneqq
    \left\{\rho \in \L1 (\reali; \reali^n) \colon \norma{\rho}_{\L1 (\reali; \reali^n)}  \leq M\right\} \,.
  \end{equation}
  Then, the map
  \begin{equation}
    \label{eq:6}
    \begin{array}{ccccc}
      \tilde\Pi_{v, \eta}
      & \colon
      & \tilde X_M
      & \to
      & \W2\infty{(\reali;\reali^n)}
      \\
      &
      & \rho
      & \mapsto
      & v\left(\cdot, (\eta * \rho) (\cdot)\right) \,,
    \end{array}
  \end{equation}
  where we denote
  \begin{equation}
    \label{eq:7}
    v\left(x, (\eta*\rho) (x)\right)
    =
    \left(
      v_i
      \left(x,
        (\eta_{i1} * \rho_1 ) (x),
        (\eta_{i2} * \rho_2 ) (x),
        \ldots,
        (\eta_{in} * \rho_n ) (x)
      \right)
    \right)_{i=1, \ldots,n},
  \end{equation}
  is well defined and Lipschitz continuous with
  \begin{equation}
    \label{eq:16}
    \!\!\Lip(\tilde{\Pi}_{v, \eta})
    \coloneqq
    \left( 7 + (6+M) M \norma{\eta}_{\W2\infty(\reali;\reali^{n \times n})} \right)
    \, \norma{v}_{\mathcal{V}^n}
    \, \norma{\eta}_{\W2\infty(\reali;\reali^{n \times n})}.
  \end{equation}
  Moreover, $\tilde\Pi_{v, \eta}$ is uniformly bounded by the quantity
  \begin{align}
    \label{stima_Q}
    Q_v
    & \coloneqq
      \big( 3 + M\norma{\eta}_{\L\infty(\reali;\reali^{n \times n})} + 3M \norma{\partial_x\eta}_{\L\infty(\reali;\reali^{n \times n})}
    \\
    \nonumber
    &\qquad +
      M^2 \norma{\partial_x\eta}^2_{\L\infty(\reali;\reali^{n \times n})}
      + M \norma{\partial_{xx}^2\eta}_{\L\infty(\reali;\reali^{n \times n})}
      \big) \norma{v}_{\mathcal{V}^n} \,.
  \end{align}
\end{lemma}

\begin{proof}
  To prove that $\tilde\Pi_{v, \eta}$ is well defined, let
  $\rho \in \L1(\reali; \reali^n )$, $i=1, \ldots, n$, and compute:
  \begin{eqnarray}
    \nonumber
    &
    & \norma{(\tilde{\Pi}_{v, \eta}\rho)_i}_{\L\infty(\reali;\reali)}
    \\
    \nonumber
    & \leq
    & \norma{v_i\left(\cdot, \eta_i * \rho(\cdot)\right)- v_i(\cdot, 0)}_{\L\infty(\reali;\reali)}
      +
      \norma{v_i(\cdot, 0)}_{\L\infty (\reali; \reali)}
    \\
    \nonumber
    & \leq
    & \norma{\nabla_\rho v_i}_{\L\infty (\reali\times\reali^n; \reali^n)} \;
      \norma{\eta_i}_{\L\infty(\reali;\reali^n)} \;
      \norma{\rho}_{\L1(\reali;\reali^n)}
      + \norma{v_i(\cdot, 0)}_{\L\infty (\reali; \reali)}
    \\
    \label{eq:11}
    & \leq
    & \left(1+M\norma{\eta_i}_{\L\infty (\reali;\reali^n)}\right)
      \norma{v_i}_{\mathcal{V}}\,;
    \\
    \nonumber
    &
    & \norma{\partial_x(\tilde{\Pi}_{v, \eta}\rho)_i}_{\L\infty(\reali;\reali)}
    \\
    \nonumber
    & \leq
    & \norma{\partial_x v_i}_{\L\infty (\reali\times\reali^n; \reali)}
      +
      \norma{
      \nabla_\rho v_i\left(\cdot, \eta_i*\rho (\cdot)\right) \;
      (\partial_x\eta_i * \rho) (\cdot)
      }_{\L\infty(\reali; \reali)}
    \\
    \nonumber
    & \leq
    & \norma{\partial_x v_i}_{\L\infty (\reali\times\reali^n; \reali)}
      +
      \norma{\nabla_\rho v_i}_{\L\infty(\reali\times\reali^n;\reali^n)} \;
      \norma{\partial_x\eta_i}_{\L\infty(\reali;\reali^n)} \;
      \norma{\rho}_{\L1(\reali;\reali^n)}
    \\
    \label{eq:12}
    & \leq
    & \left(1+M\norma{\partial_x\eta_i}_{\L\infty (\reali;\reali^n)}\right)
      \norma{v_i}_{\mathcal{V}}\,;
    \\ 
    \nonumber
    &
    & \norma{\partial_{xx}^2 (\tilde{\Pi}_{v, \eta}\rho)_i
      }_{\L\infty(\reali;\reali)}
    \\
    \nonumber
    & \leq
    & \norma{\partial^2_{xx} v_i}_{\L\infty (\reali\times \reali^n; \reali)}
      + 2 \norma{\partial_x \nabla_\rho v_i\left(\cdot, \eta_i*\rho (\cdot)\right) \; (\partial_x \eta_i *\rho) (\cdot)}_{\L\infty (\reali; \reali)}
    \\
    \nonumber
    &
    & +
      \norma{D^2_{\rho\rho} v_i\left(\cdot, \eta_i*\rho (\cdot)\right) \; \left(\partial_x\eta_i*\rho (\cdot)\right)^2
      +
      \nabla_\rho v_i\left(\cdot, \eta_i*\rho (\cdot)\right) \; (\partial_{xx}^2\eta_i*\rho) (\cdot)
      }_{\L\infty(\reali;\reali)}
    \\
    \nonumber
    & \leq
    & \norma{\partial^2_{xx} v_i}_{\L\infty (\reali\times \reali^n; \reali)}
      + 2
      \norma{\partial_x \nabla_\rho v_i}_{\L\infty (\reali\times\reali^n; \reali^n)} \; \norma{\partial_x \eta_i}_{\L\infty (\reali;\reali^n)}
      \norma{\rho}_{\L1 (\reali; \reali^n)}
    \\
    \nonumber
    &
    & +
      \norma{D^2_{\rho\rho} v_i}_{\L\infty(\reali\times\reali^n; \reali^{n\times n})} \;
      \norma{\partial_x\eta_i}_{\L\infty(\reali;\reali^n)}^2 \;
      \norma{\rho}^2_{\L1(\reali;\reali^n)}
    \\
    \nonumber
    &
    & \quad
      +
      \norma{\nabla_\rho v_i}_{\L\infty(\reali\times\reali^n; \reali^n)} \;
      \norma{\partial_{xx}^2 \eta_i}_{\L\infty(\reali;\reali^n)} \;
      \norma{\rho}_{\L1(\reali;\reali^n)}
    \\
    \label{eq:13}
    & \leq
    & \left(1
      + 2M \norma{\partial_x\eta_i}_{\L\infty (\reali;\reali^n)}
      + M^2 \norma{\partial_x\eta_i}_{\L\infty (\reali;\reali^n)}^2
      + M \norma{\partial^2_{xx}\eta_i}_{\L\infty (\reali;\reali^n)}
      \right) \norma{v_i}_{\mathcal{V}}\,;
  \end{eqnarray}
  hence, $\tilde\Pi_{v, \eta}$ is well defined.

  To prove Lipschitz continuity, let
  $r, \rho \in \L1(\reali; \reali^n )$, $i=1, \ldots, n$ and evaluate
  \begin{eqnarray*}
    &
    & \norma{(\tilde{\Pi}_{v, \eta}r)_i - (\tilde{\Pi}_{v, \eta}\rho)_i}_{\L\infty(\reali;\reali)}
    \\
    & \leq
    & \norma{\nabla_\rho v_i}_{\L\infty(\reali\times\reali^n;\reali^n)} \;
      \norma{\eta_i *(r-\rho)}_{\L\infty(\reali; \reali^n)}
    \\
    & \leq
    & \norma{v_i}_{\mathcal{V}} \;
      \norma{\eta_i}_{\L\infty(\reali;\reali^n)} \;
      \norma{r-\rho}_{\L1(\reali; \reali^n)} \,;
    \\
    &
    & \norma{\partial_x (\tilde{\Pi}_{v, \eta}r)_i-\partial_x (\tilde{\Pi}_{v, \eta}\rho)_i
      }_{\L\infty(\reali; \reali)}
    \\
    & \leq
    &
      \norma{\partial_x v_i\left(\cdot, (\eta_i * r) (\cdot)\right)
      - \partial_x v_i\left(\cdot, (\eta_i * \rho) (\cdot)\right)
      }_{\L\infty (\reali; \reali)}
    \\
    &
    & + \norma{
      \nabla_\rho v_i\left(\cdot, (\eta_i * r) (\cdot)\right) \; (\partial_x \eta_i * r) (\cdot)
      -
      \nabla_\rho v_i\left(\cdot,(\eta_i * \rho) (\cdot)\right) \; (\partial_x\eta_i *\rho) (\cdot)
      }_{\L\infty(\reali; \reali)}
    \\
    & \leq
    &
      \norma{\partial_x \nabla_\rho v_i}_{\L\infty (\reali\times\reali^n;\reali^n)}
      \norma{\eta_i}_{\L\infty(\reali;\reali^n)} \norma{ r-\rho}_{\L1 (\reali;\reali^n)}
    \\
    &
    & + \norma{ \nabla_\rho v_i}_{\L\infty(\reali \times \reali^n; \reali^n)} \norma{\partial_x \eta_i}_{\L\infty(\reali;\reali^n)} \norma{r-\rho}_{\L1(\reali;\reali^n)}
    \\
    &
    & + M \norma{D^2_{\rho \rho} v_i}_{\L\infty(\reali \times \reali^n; \reali^{n \times n})} \norma{\eta_i}_{\L\infty(\reali;\reali^n)} \norma{\partial_x \eta_i}_{\L\infty(\reali;\reali^n)} \norma{r- \rho}_{\L1(\reali;\reali^n)}
    \\
    & \leq
    & \left( 2 + M \norma{\eta_i}_{\W2\infty(\reali;\reali^n)}\right) \norma{v_i}_{\mathcal{V}}\norma{\eta_i}_{\W2\infty(\reali;\reali^n)}\norma{r- \rho}_{\L1(\reali;\reali^n)} \,;
    \\
    &
    & \norma{ \partial_{xx}^2 \left(\tilde{\Pi}_{v, \eta}r\right) _i-\partial_{xx}^2 \left(\tilde{\Pi}_{v, \eta}\rho\right) _i }_{\L\infty(\reali;\reali)}
    \\
    & \leq
    &
      \norma{\partial_{xx}^2 v_i(\cdot, (\eta_i * r)(\cdot)) - \partial_{xx}^2 v_i(\cdot, (\eta_i * \rho)(\cdot)) }_{\L\infty(\reali;\reali)}
    \\
    &
    & + 2\norma{\partial_x \nabla_\rho v_i(\cdot, (\eta_i * r)(\cdot))\left( (\partial_x \eta_i * (r -\rho)(\cdot)\right)}_{\L\infty(\reali;\reali)}
    \\
    &
    & + 2\norma{\left( \partial_x \nabla_\rho v_i(\cdot, (\eta_i * r)(\cdot)) - \partial_x \nabla_\rho v_i(\cdot, (\eta_i * \rho)(\cdot))\right) (\partial_x \eta_i * \rho)(\cdot)}_{\L\infty(\reali;\reali)}
    \\
    &
    & + \norma{D_{\rho \rho}^2 v_i(\cdot, (\eta_i * r)(\cdot))\left( (\partial_x \eta_i * r(\cdot))^2 - (\partial_x \eta_i * \rho(\cdot))^2 \right)}_{\L\infty(\reali;\reali)}
    \\
    &
    &+ \norma{\left( D_{\rho \rho}^2 v_i(\cdot, (\eta_i * r)(\cdot)) - D_{\rho \rho}^2 v_i(\cdot, (\eta_i * \rho)(\cdot))\right)(\partial_x \eta_i * \rho(\cdot))^2}_{\L\infty(\reali;\reali)}
    \\
    &
    &+ \norma{\nabla_\rho v_i(\cdot, (\eta_i * r(\cdot)))\left( \partial_{xx}^2 \eta_i * (r-\rho)(\cdot) \right) }_{\L\infty(\reali;\reali)}
    \\
    &
    &+ \norma{\left( \nabla_\rho v_i(\cdot, (\eta_i * r(\cdot))) - \nabla_\rho v_i(\cdot, (\eta_i * \rho(\cdot)))\right)(\partial_{xx}^2 \eta_i * \rho(\cdot))  }_{\L\infty(\reali;\reali)}
    \\
    & \leq
    & \Big( \norma{\partial_{xx}^2 \nabla_\rho v_i}_{\L\infty(\reali \times \reali^n; \reali^n)} \norma{\eta_i}_{\L\infty(\reali;\reali^n)}
    \\
    &
    &+ 2\norma{\partial_x \nabla_\rho v_i}_{\L\infty(\reali \times \reali^n; \reali^n)}\norma{\partial_x \eta_i}_{\L\infty(\reali;\reali^n)}
    \\
    &
    &+ 2 \norma{\partial_x D_{\rho \rho}^2 v_i}_{\L\infty(\reali \times \reali^n; \reali^{n\times n})} \norma{\eta_i}_{\L\infty(\reali;\reali^n)}\norma{\partial_x \eta_i}_{\L\infty(\reali;\reali^n)}M
    \\
    &
    &+ \norma{D_{\rho \rho}^2 v_i}_{\L\infty(\reali \times \reali^n; \reali^{n\times n})} \norma{\partial_x \eta_i}_{\L\infty(\reali;\reali^n)}^22M
    \\
    &
    &+ \norma{D^3_{\rho \rho} v_i}_{\L\infty(\reali \times \reali^n; \reali^{n\times n \times n })}\norma{\eta_i}_{\L\infty(\reali;\reali^n)} \norma{\partial_x \eta_i}^2_{\L\infty(\reali;\reali^n)} M^2
    \\
    &
    &+ \norma{\nabla_\rho v_i}_{\L\infty(\reali \times \reali^n; \reali^{n})} \norma{\partial_{xx}^2 \eta_i}_{\L\infty(\reali;\reali^n)}
    \\
    &
    &+ \norma{D^2_{\rho \rho}v_i}_{\L\infty(\reali \times \reali^n; \reali^{n\times n})} \norma{\eta_i}_{\L\infty(\reali;\reali^n)} \norma{\partial^2_{xx}\eta_i}_{\L\infty(\reali;\reali^n)}M \Big)
    \\
    &
    & \times \norma{r-\rho}_{\L1(\reali;\reali^n)}
    \\
    & \leq
    & \left( 4 + 5M\norma{\eta_i}_{\W2\infty(\reali;\reali^n)}+M^2 \norma{\eta_i}_{\W2\infty(\reali;\reali^n)} \right) \norma{v_i}_{\mathcal{V}} \norma{\eta_i}_{\W2\infty(\reali;\reali^n)} \norma{r-\rho}_{\L1(\reali;\reali^n)}.
  \end{eqnarray*}
  Thus, $\tilde{\Pi}_{v, \eta}$ is Lipschitz continuous with a
  Lipschitz constant depending on $v,M$ and $\eta$.

  The uniform bound of $\tilde\Pi_{v,\eta}$ directly follows
  from~\eqref{eq:11}, \eqref{eq:12} and~\eqref{eq:13}.
\end{proof}

\noindent We now investigate the dependence of the map
$\tilde{\Pi}_{v, \eta}$ on $v$ and $\eta$.

\begin{proposition}
  \label{prop:5}
  Fix $\rho \in \tilde X_M$ as defined in~\eqref{eq:2}. With the
  notation~\eqref{eq:6}--\eqref{eq:7}, define
  \begin{equation}
    \label{eq:8}
    \begin{array}{c@{\,}c@{\,}lcc}
      \mathcal{V}^n
      & \times
      & \W2\infty (\reali; \reali^{n\times n})
      & \to
      & \W2\infty (\reali; \reali^n)
      \\
      w
      & ,
      & \;\xi
      & \mapsto
      & \tilde{\Pi}_{w, \xi} \rho \,.
    \end{array}
  \end{equation}
  Then,
  \begin{enumerate}[label=\textbf{(\arabic*)}]
  \item \label{item:10} For all
    $\xi \in \W2\infty (\reali; \reali^n)$, the map
    $v \mapsto \tilde{\Pi}_{v, \xi} \rho$ is Lipschitz continuous:
    $\forall v,w \in \mathcal{V}^n$
    \begin{equation*}
      \norma{\tilde{\Pi}_{v, \xi}\rho - \tilde{\Pi}_{w, \xi}\rho }_{\W2\infty(\reali;\reali^n)} \leq
      \left(3+5M\norma{\xi}_{\W2\infty(\reali;\reali^n)} + M^2\norma{\xi}^2_{\W2\infty(\reali;\reali^n)}\right) \norma{v-w}_{\mathcal{V}^n}.
    \end{equation*}
  \item \label{item:11} For all $v \in \mathcal{V}^n$, the map
    $\xi \mapsto \tilde{\Pi}_{v, \xi} \rho$ is locally Lipschitz
    continuous: $\forall \xi, \eta \in \W2\infty(\reali;\reali^n)$
    \begin{align*}
      &\norma{\tilde{\Pi}_{v, \xi}\rho - \tilde{\Pi}_{v, \eta}\rho }_{\W2\infty(\reali;\reali^n)}
      \\
      \leq
      &
        \left(7 + 4M\norma{\eta}_{\W2\infty(\reali;\reali^n)} + 3M\norma{\xi}_{\W2\infty(\reali;\reali^n)} + (M+M^2)\norma{\xi}^2_{\W2\infty(\reali;\reali^n)}\right)
      \\
      &\times
        M\norma{v}_{\mathcal{V}^n}\norma{\xi-\eta}_{\W2\infty(\reali;\reali^n)}.
    \end{align*}
  \end{enumerate}
\end{proposition}

\begin{proof}
  Fix $\xi \in \W2\infty (\reali;\reali^{n\times n})$,
  $\rho \in \tilde{X}_M$ and $v, w \in \mathcal{V}$.  For
  $i=1, \ldots, n$, one obtains
  \begin{eqnarray*}
    &
    & \norma{(\tilde{\Pi}_{v, \xi} \rho)_i - (\tilde{\Pi}_{w, \xi} \rho)_i}_{\L\infty(\reali;\reali)}
    \\
    & \leq
    & \norma{v_i (\cdot, 0) - w_i (\cdot, 0)}_{\L\infty (\reali;\reali)}
      +
      \norma{\nabla_\rho v_i - \nabla_\rho w_i}_{\L\infty (\reali\times\reali^n; \reali^n)}
      \norma{\xi_i}_{\L\infty (\reali; \reali^n)}
      \norma{\rho}_{\L1 (\reali; \reali^n)}
    \\
    & \leq
    & \left(1 + M \, \norma{\xi_i}_{\L\infty (\reali; \reali^n)}\right)
      \norma{v_i - w_i}_{\mathcal{V}} \, .
  \end{eqnarray*}
  Moreover, for $i=1, \ldots, n$,
  \begin{eqnarray*}
    &
    & \norma{\partial_x (\tilde{\Pi}_{v, \xi} \rho)_i -
      \partial_x(\tilde{\Pi}_{w, \xi} \rho)_i}_{\L\infty(\reali;\reali)}
    \\
    & \leq
    & \norma{\partial_x v_i\left(\cdot, \xi_i *\rho(\cdot)\right) -
      \partial_x w_i\left(\cdot, \xi_i *\rho(\cdot)\right) }_{\L\infty(\reali;\reali)}
    \\
    &
    & +
      \norma{\left( \nabla_\rho{v_i}\left(\cdot, \xi_i*\rho (\cdot)\right) - \nabla_\rho{w_i}\left(\cdot, \xi_i*\rho (\cdot)\right)\right)\left(\partial_x \xi_i *\rho(\cdot)\right)}_{\L\infty(\reali;\reali)}
    \\
    & \leq
    &
      \norma{\partial_x v_i - \partial_x w_i}_{\L\infty(\reali \times \reali^n; \reali)}
      +
      \norma{\rho}_{\L1(\reali;\reali^n)} \norma{\partial_x \xi_i}_{\L\infty(\reali;\reali^n)}\norma{\nabla_\rho v_i -\nabla_\rho w_i}_{\L\infty(\reali \times \reali^n,\reali^n)}
    \\
    & \leq
    & \left(1 + M \, \norma{\partial_x \xi_i}_{\L\infty (\reali; \reali^n)}\right)
      \norma{v_i - w_i}_{\mathcal{V}}
  \end{eqnarray*}
  and
  \begin{eqnarray*}
    &
    & \norma{\partial^2_{xx}(\tilde{\Pi}_{v, \xi} \rho)_i - \partial^2_{xx} (\tilde{\Pi}_{w,\xi} \rho)_i}_{\L\infty(\reali;\reali)}
    \\
    & =
    &
      \norma{\partial_{xx}^2v_i(\cdot, (\xi_i*\rho)(\cdot)) - \partial_{xx}^2w_i(\cdot, (\xi_i*\rho)(\cdot))}_{\L\infty(\reali ; \reali)}
    \\
    &
    &+
      \norma{\left( \partial_x\nabla_\rho v_i(\cdot, (\xi_i*\rho)(\cdot))  -\partial_x\nabla_\rho w_i(\cdot, (\xi_i*\rho)(\cdot)) \right) 2(\partial_x \xi_i *\rho)(\cdot) }_{\L\infty(\reali ; \reali)}
    \\
    &
    &+
      \norma{\left(\nabla_{\rho\rho}^2 v_i(\cdot, (\xi_i* \rho)(\cdot)) - \nabla_{\rho\rho}^2 w_i(\cdot, (\xi_i* \rho)(\cdot))\right) (\partial_x \xi_i*\rho)^2 (\cdot)}_{\L\infty(\reali ; \reali)}
    \\
    &
    &+
      \norma{\left(\nabla_\rho v_i(\cdot, (\xi_i*\rho)(\cdot)) - \nabla_\rho w_i(\cdot, (\xi_i*\rho)(\cdot))\right) (\partial_{xx}^2 \xi_i * \rho)(\cdot)}_{\L\infty(\reali ; \reali)}
    \\
    & \leq
    &
      \norma{\partial_{xx}^2 v_i - \partial_{xx}^2 w_i}_{\L\infty(\reali\times\reali^n;\reali )}
    \\
    &
    &+
      \norma{\partial_x \nabla_\rho v_i - \partial_x \nabla_\rho w_i }_{\L\infty(\reali\times\reali^n;\reali^n )}2\norma{\partial_x \xi_i}_{\L\infty(\reali;\reali^n )}M
    \\
    &
    &+
      \norma{\nabla_{\rho\rho}^2 v_i - \nabla_{\rho\rho}^2 w_i}_{\L\infty(\reali\times\reali^n;\reali^{n \times n} )}\norma{\partial_x \xi_i}_{\L\infty(\reali;\reali^n )}^2 M^2
    \\
    &
    &+
      \norma{\nabla_\rho v_i - \nabla_\rho w_i}_{\L\infty(\reali\times\reali^n;\reali^n )}\norma{\partial_{xx}^2 \xi_i }_{\L\infty(\reali;\reali^n )} M
    \\
    & \leq
    &\left(
      1
      +2 M \norma{\partial_x \xi_i}_{\L\infty (\reali; \reali^n)}
      + M^2 \norma{\partial_x \xi_i}_{\L\infty (\reali; \reali^n)}^2
      + M \norma{\partial^2_{xx} \xi_i}_{\L\infty (\reali; \reali^n)}
      \right)
      \norma{v_i - w_i}_{\mathcal{V}} \,,
  \end{eqnarray*}
  completing the proof of~\ref{item:10}.

  Now, we prove that the function
  $\xi \mapsto \tilde{\Pi}_{v, \xi}\rho $ is locally Lipschitz
  continuous.  Fix $v \in \mathcal{V}^n$, $\rho \in \tilde{X}_M$ and
  the functions $\eta, \xi \in \W2\infty(\reali;\reali^{n\times n})$.
  It is immediate to verify that for $i=1, \ldots, n$,
  \begin{eqnarray*}
    \norma{(\tilde{\Pi}_{v,\eta}\rho)_i - (\tilde{\Pi}_{v,\xi} \rho)_i}_{\L\infty(\reali;\reali)}
    & =
    &
      \norma{ v_i(\cdot, \eta_i * \rho(\cdot)) - v_i(\cdot, \xi_i * \rho(\cdot))}_{\L\infty(\reali;\reali)}
    \\
    & \leq
    & M \norma{\nabla_\rho v_i}_{\L\infty(\reali \times \reali^n;\reali^n)} \norma{\eta_i - \xi_i}_{\L\infty(\reali;\reali^n)}
  \end{eqnarray*}
  and
  \begin{eqnarray*}
    &
    &\norma{\partial_x(\tilde{\Pi}_{v,\eta} \rho)_i - \partial_x(\tilde{\Pi}_{v,\xi} \rho)_i}_{\L\infty(\reali;\reali)}
    \\
    & \leq
    &
      \norma{\partial_x v_i(\cdot, \eta_i * \rho(\cdot)) - \partial_x v_i(\cdot, \xi_i * \rho(\cdot))}_{\L\infty(\reali;\reali)}
    \\
    &
    &+
      \norma{\left(\nabla_\rho v_i(\cdot, \eta_i *\rho(\cdot)) - \nabla_\rho v_i(\cdot, \xi_i *\rho(\cdot))\right) (\partial_x \eta_i *\rho (\cdot))}_{\L\infty(\reali;\reali)}
    \\
    &
    &+
      \norma{\nabla_\rho v_i(\cdot, \xi_i *\rho(\cdot)) \left( (\partial_x \eta_i - \partial_x \xi_i)*\rho (\cdot)\right) }_{\L\infty(\reali;\reali)}
    \\
    & \leq
    &
      \norma{\partial_x \nabla_\rho v_i}_{\L\infty(\reali \times \reali^n; \reali^n)} \norma{\eta_i - \xi_i}_{\L\infty(\reali;\reali^n)}M
    \\
    &
    &+
      \norma{D_{\rho\rho}^2 v_i}_{\L\infty(\reali \times \reali^n; \reali^{n\times n})} M^2 \norma{\eta_i- \xi_i}_{\L\infty(\reali;\reali^n)}\norma{\partial_x \eta_i }_{\L\infty(\reali,\reali^n)}
    \\
    &
    &+
      \norma{\nabla_\rho v_i}_{\L\infty(\reali \times \reali^n; \reali^n)} \norma{\partial_x \eta_i - \partial_x \xi_i}_{\L\infty(\reali;\reali^n)}M
    \\
    & \leq
    & M \norma{v_i}_{\mathcal{V}}
      \left(2 + M \norma{\partial_x \eta_i}_{\L\infty (\reali; \reali^n)}\right)
      \norma{\eta_i - \xi_i}_{\W1\infty (\reali; \reali^n)} \,.
  \end{eqnarray*}
  Moreover, for $i=1, \ldots, n$,
  \begin{eqnarray*}
    &
    & \norma{\partial_{xx}^2(\tilde{\Pi}_{v,\eta}  \rho)_i - \partial_{xx}^2(\tilde{\Pi}_{v,\xi}  \rho)_i}_{\L\infty(\reali;\reali)}
    \\
    & \leq
    &
      \norma{\partial_{xx}^2 v_i(\cdot, \eta_i*\rho(\cdot)) - \partial_{xx}^2 v_i(\cdot, \xi_i*\rho(\cdot))}_{\L\infty(\reali;\reali)}
    \\
    &
    &+
      \norma{2 \partial_x \nabla_\rho v_i(\cdot, \eta_i *\rho(\cdot))\left( (\partial_x \eta_i - \partial_x \xi_i)*\rho(\cdot) \right)}_{\L\infty(\reali;\reali)}
    \\
    &
    &+
      \norma{2 \left(\partial_x\nabla_\rho v_i(\cdot, \eta_i*\rho(\cdot)) - \partial_x\nabla_\rho v_i(\cdot, \xi_i*\rho(\cdot)) \right)(\partial_x \xi_i* \rho(\cdot))}_{\L\infty(\reali;\reali)}
    \\
    &
    &+
      \norma{\nabla_{\rho\rho}^2 v_i(\cdot, \eta_i*\rho(\cdot)) \left((\partial_x \eta_i *\rho(\cdot))^2 - (\partial_x \xi_i *\rho(\cdot))^2 \right)}_{\L\infty(\reali;\reali)}
    \\
    &
    &+
      \norma{\left( \nabla_{\rho\rho}^2v_i(\cdot, \eta_i*\rho(\cdot)) - \nabla_{\rho\rho}^2v_i(\cdot, \xi_i*\rho(\cdot))
      \right)(\partial_x \xi_i*\rho(\cdot))^2}_{\L\infty(\reali;\reali)}
    \\
    &
    &+
      \norma{\nabla_\rho v_i(\cdot, \eta_i*\rho(\cdot))\left( (\partial_{xx}^2 \eta_i - \partial_{xx}^2 \xi_i)*\rho(\cdot)\right)}_{\L\infty(\reali;\reali)}
    \\
    &
    &+
      \norma{\left( \nabla_\rho v_i(\cdot, \eta_i*\rho(\cdot)) - \nabla_\rho v_i(\cdot, \xi_i*\rho(\cdot))\right)\left(\partial_{xx}^2\xi_i*\rho(\cdot)\right)}_{\L\infty(\reali;\reali)}
    \\
    & \leq
    &
      \norma{\partial_{xx}^2 \nabla_\rho v_i}_{\L\infty(\reali \times \reali^n; \reali^n)} \norma{\eta_i -\xi_i}_{\L\infty(\reali;\reali^n)} M
    \\
    &
    &+
      2\norma{\partial_x \nabla_\rho v_i}_{\L\infty(\reali \times \reali^n; \reali^n)}\norma{\partial_x \eta_i - \partial_x \xi_i}_{\L\infty(\reali;\reali^n)}M
    \\
    &
    &+
      2\norma{\partial_x D_{\rho\rho}^2 v_i}_{\L\infty(\reali \times \reali^n; \reali^{n \times n})}\norma{\eta_i - \xi_i}_{\L\infty(\reali;\reali^n)}M^2 \norma{\partial_x \xi_i}_{\L\infty(\reali;\reali^n)}
    \\
    &
    &+
      \norma{D^2_{\rho \rho} v_i}_{\L\infty(\reali \times \reali^n; \reali^{n \times n})} \norma{\partial_x \eta_i - \partial_x \xi_i}_{\L\infty(\reali;\reali^n)} M^2 \left( \norma{\partial_x \eta_i}_{\L\infty(\reali;\reali^n)} + \norma{\partial_x \xi_i}_{\L\infty(\reali;\reali^n)}\right)
    \\
    &
    &+
      \norma{D_{\rho\rho\rho}^3 v_i}_{\L\infty(\reali \times \reali^n; \reali^{n \times n \times n})}\norma{\eta_i - \xi_i}_{\L\infty(\reali;\reali^n)}\norma{\partial_x \xi_i}_{\L\infty(\reali;\reali^n)}^2 M^3
    \\
    &
    &+
      \norma{\nabla_\rho v_i}_{\L\infty(\reali \times \reali^n; \reali^n)} \norma{\partial_{xx}^2 \eta_i- \partial_{xx}^2 \xi_i}_{\L\infty(\reali;\reali^n)} M
    \\
    &
    &+
      \norma{D_{\rho\rho}^2 v_i}_{\L\infty(\reali \times \reali^n; \reali^{n \times n})}\norma{\eta_i -\xi_i}_{\L\infty(\reali;\reali^n)}\norma{\partial_{xx}^2 \xi_i}_{\L\infty(\reali;\reali^n)}M^2
    \\
    & \leq
    & \left(
      4 + 3 M
      \left(\norma{\partial_x \eta_i}_{\L\infty (\reali;\reali^n)}
      {+}
      \norma{\partial_x \xi_i}_{\L\infty (\reali;\reali^n)}\right)
      {+} M^2 \norma{\partial_x\xi_i}_{\L\infty (\reali;\reali^n)}^2
      {+} M \norma{\partial^2_{xx}\xi_i}_{\L\infty (\reali;\reali^n)}^2
      \right)
    \\
    &
    & \times M \norma{v_i}_{\mathcal{V}}
      \norma{\eta_i-\xi_i}_{\W2\infty (\reali; \reali^n)} \,,
  \end{eqnarray*}
  proving~\ref{item:11}.
\end{proof}

In the next step we introduce the time dependent version of the map
$\tilde{\Pi}_{v,\eta}$ defined by~\eqref{eq:6}--\eqref{eq:7} on the
time interval $[0,T]$, for a fixed $T > 0$. To this aim, recall
$\tilde X_M$ as in~\eqref{eq:2} and introduce
\begin{equation}
  \label{eq:10}
  X_M
  \coloneqq
  \C0([0,T]; \tilde{X}_M)
  \quad \mbox{ and } \quad
  \norma{\rho}_{X_M}
  \coloneqq
  \sup_{t \in [0,T]} \norma{\rho (t)}_{\L1 (\reali;\reali^n)} \,.
\end{equation}

\begin{lemma}
  \label{lem:Pi}
  Let $v \in \C0 ([0,T]; \mathcal{V}^n)$ and~\ref{item:2} hold.  Then,
  the map
  \begin{displaymath}
    \begin{array}{ccccl}
      \Pi_{v, \eta}
      & \colon
      & X_M
      & \to
      & \C0\left([0,T]; \W2\infty{(\reali;\reali^n)}\right)
      \\
      \\
      &
      & \rho
      & \mapsto
      & \left[ t \mapsto \tilde\Pi_{v (t), \eta}\left( \rho (t)\right)\right]
        \qquad \mbox{i.e. }
        \left((\Pi_{v,\eta}\rho) (t)\right) (x)
        =
        v\left(t, x, \left(\eta * \rho (t)\right) (x)\right)
    \end{array}
  \end{displaymath}
  is well defined, bounded and satisfies for $r,\rho \in X_M$ the
  Lipschitz estimates
  \begin{eqnarray}
    \label{eq:lip-1}
    \norma{\Pi_{v, \eta} r - \Pi_{v, \eta} \rho}_{\C0([0,T]; \W2\infty(\reali; \reali^n)) }
    & \leq
    & \Lip (\Pi_{v, \eta}) \; \norma{r-\rho}_{X_M}
    \\
    \label{eq:lip-2}
    \norma{\Pi_{v, \eta} r - \Pi_{v, \eta} \rho}_{\L1([0,T]; \W2\infty(\reali; \reali^n)) }
    & \leq
    & \Lip (\Pi_{v, \eta}) \;
      \norma{r-\rho}_{\L1 ([0,T]\times\reali;\reali^n)}
  \end{eqnarray}
  where $\Lip (\Pi_{v, \eta})=\sup_{t \in [0,T]} \Lip \left ( \tilde \Pi_{v(t), \eta} \right)$.
\end{lemma}

\begin{proof}
  Fix $\rho \in X_M$, $\bar{t} \in [0,T]$ and consider a sequence $\{ t_n\}_n \subset [0,T]$ converging to $\bar t$. Taking advantage of the
  Lipschitz regularity of $\tilde{\Pi}_{v,\eta}$, by
  Lemma~\ref{lem:Pitilde} and~\ref{item:10} in
  Proposition~\ref{prop:5}, one can get
  \begin{eqnarray*}
    &
    &\norma{ \tilde{\Pi}_{v (t_n),\eta}(\rho(t_n)) -
      \tilde{\Pi}_{v (\bar t),\eta}(\rho(\bar t)) }_{\W2\infty(\reali;\reali^n)}
    \\
    & \leq
    & \norma{ \tilde{\Pi}_{v (t_n),\eta}(\rho(t_n)) -
      \tilde{\Pi}_{v (\bar t),\eta}(\rho(t_n)) }_{\W2\infty(\reali;\reali^n)}
      + \norma{ \tilde{\Pi}_{v (\bar t),\eta}(\rho(t_n)) -
      \tilde{\Pi}_{v (\bar t),\eta}(\rho(\bar t)) }_{\W2\infty(\reali;\reali^n)}
    \\
    &\leq
    & C (M, \norma{\eta}_{\W2\infty (\reali;\reali^{n\times n})}) \,
      \norma{v (t_n) - v (\bar t)}_{\mathcal{V}^n}
      +
      \Lip{\tilde{\Pi}_{v (\bar t), \eta}} \,
      \norma{ \rho(t_n) - \rho(\bar t)}_{\L1(\reali;\reali^n)}
  \end{eqnarray*}
  which tends to zero as $n \to +\infty$, due to the
  continuity in time of $v$ and $\rho$. Hence, $\Pi_{v, \eta}$ is well defined.

  The boundedness of $\Pi_{v, \eta}$ follows from the same property of
  $v$ by~\ref{item:1} and of $\tilde\Pi_{v (t), \eta}$ as proved in
  Lemma~\ref{lem:Pitilde} for all $t \in [0,T]$.

  To prove the Lipschitz estimates on $\Pi_{v, \eta}$, let
  $r,\rho \in X_M$ and evaluate
  \begin{eqnarray*}
    \norma{\Pi_{v, \eta} r - \Pi_{v, \eta} \rho}_{\C0([0,T]; \W2\infty(\reali; \reali^n)) }
    & =
    & \sup_{t \in [0,T]} \norma{\tilde\Pi_{v (t), \eta} \left(r (t)\right)
      - \tilde\Pi_{v (t), \eta} \left(\rho (t)\right)}_{\W2\infty (\reali; \reali^n)}
    \\
    & \leq
    & \sup_{t \in [0,T]} \Lip (\tilde\Pi_{v (t), \eta}) \norma{r (t) - \rho (t)}_{\L1 (\reali;\reali^n)}
    \\
    & \leq
    & \Lip (\Pi_{v, \eta}) \; \norma{r-\rho}_{X_M} \,,
  \end{eqnarray*}
  where
  $\Lip (\Pi_{v, \eta}) \coloneqq \sup_{t \in [0,T]} \Lip (\tilde\Pi_{v (t),
    \eta})$ is finite by~\eqref{eq:16} and $v \in \C0([0, T]; \mathcal{V}^n)$,
    proving~\eqref{eq:lip-1}. An analogous procedure gives the other
  estimate~\eqref{eq:lip-2}.
\end{proof}

For later use, recall that, given
$w \in \C0([0,T];\W1\infty(\reali;\reali^n))$ and
$\rho_o \in (\L1 \cap \BV) (\reali; \reali^n)$, the Cauchy problems
\begin{equation}\label{eq:21}
  \begin{cases}
    \partial_t \rho_i + \partial_x (\rho_i\, w_i)=0\\
    \rho_i(0,x)=(\rho_o)_i(x)
  \end{cases}
  \qquad i=1,\ldots, n,
\end{equation}
admit the (Lagrangian, \rev{see~\cite[\S~2.5]{MR3057143}}) solutions $\rho_1, \ldots, \rho_n$ such that
\begin{equation}
  \label{eq:4}
  \rho_i(t,x)= (\rho_o)_i(X_i(0;t,x)) \; \exp \left( -\int_{0}^t \partial_x w_i(s,X_i(s;t,x)) \, \d{s}\right)
\end{equation}
where the characteristic $t \mapsto X_i(t;t_o,x_o)$ is the solution to
the Cauchy problem
\begin{equation}
  \label{eq:24}
  \left \{
    \begin{array}{ll}
      \displaystyle \dot{x} = w_i(t, x) & \ t \in \reali_+ \\
      \displaystyle x(t_o)=x_o \,.
    \end{array}
  \right.
\end{equation}
The maps $\rho_1, \ldots, \rho_n$ are also Kružkov solutions
to~\eqref{eq:21} in the sense
of~\cite[Definition~1]{MR267257}. Indeed, several results in the
literature, see for instance~\cite[Lemma 5]{MR4371486}
or~\cite[Corollary~II.1]{MR1022305}, ensure that for a Cauchy problem
of the type \eqref{eq:21}, the concepts of weak and entropy (or
Kružkov) solutions coincide.

Useful relations (see~\cite[Chapter~3]{MR2347697} or~\cite[Lemma
2.6]{MR3670045}) which we exploit below are:
\begin{eqnarray}
  \nonumber
  \partial_t X_i(t;t_o,x_o)
  & =
  & w_i\left(t, X_i(t;t_o,x_o)\right)
  \\
  \nonumber
  \partial_{t_o} X_i(t;t_o,x_o)
  & =
  & - w_i(t_o,x_o) \; \exp\left ( \int_0^t \partial_x w_i(s, X_i(s;t_o,x_o)) \d{s} \right)
  \\
  \label{eq:1}
  \partial_{x_o} X_i(t;t_o,x_o)
  & =
  & \rev{\exp\left(\int_{t_o}^t
    \partial_x w_i \left(s, X_i (s;t_o,x_o)\right)\d{s}\right)}
\end{eqnarray}
where $i \in \{1, \ldots, n\}$.

\begin{lemma}
  \label{lem:Sigma}
  Let $\rho_o \in (\L1 \cap \BV) (\reali; \reali^n)$, fix ${\tilde Q} > 0$ and
  call
  \begin{equation}
    \label{eq:5}
    W_{\tilde Q}
    \coloneqq
    \left\{ w \in \C0\left([0,T]; \W2\infty(\reali; \reali^n)\right) \colon
      \norma{w}_{\C0([0,T]; \W2\infty(\reali; \reali^n))}\leq {\tilde Q}\right\}.
  \end{equation}
  Define the map $\Sigma_{\rho_o}$ as
  \begin{equation}
    \label{def:Sigma}
    \begin{array}{ccccc}
      \Sigma_{\rho_o}
      & \colon
      & W_{\tilde Q}
      & \to
      & \C0\left([0,T]; \L1(\reali; \reali^N)\right)
      \\
      &
      & w
      & \mapsto
      & \rho \,,
    \end{array}
  \end{equation}
  where $\rho$ is the solution to~\eqref{eq:21}. Then,
  \begin{enumerate}[label=\bf($\mathbf\Sigma\arabic*$)]
  \item \label{item:3} The map $\Sigma_{\rho_o}$ is well defined and
    Lipschitz continuous.
  \item \label{item:9} For all $w \in W_{\tilde Q}$,
    $t \mapsto (\Sigma_{\rho_o}w) (t)$ is locally $\L1$-Lipschitz
    continuous in time.
  \item \label{item:12} For all $w \in W_{\tilde Q}$, $t \in [0,T]$ and
    $i = 1, \ldots, n$,
    \begin{eqnarray}
      \label{eq:15}
      \norma{(\Sigma_{\rho_o} w)_i(t)}_{\L1(\reali;\reali)}
      & =
      & \norma{(\rho_o)_i}_{\L1(\reali;\reali)}
        \quad \forall t\in [0,T] ,
        \quad \forall i=1, \ldots, n \,;
      \\
      \label{eq:tv_Sigma}
      \tv\left((\Sigma_{\rho_o} w)_i(t)\right)
      & \leq
      & \left(
        \tv\left((\rho_o)_i\right)
        +
        {\tilde Q} \, t \, \norma{(\rho_o)_i}_{\L1 (\reali; \reali)}
        \right)
        e^{{\tilde Q}\, t} \,.
    \end{eqnarray}
  \item \label{item:13} For all $i=1, \ldots,n$ and for all
    $w \in W_{\tilde Q}$, if $(\rho_o)_i \geq 0$, then for all $t \in [0,T]$,
    $\rho_i (t) \geq 0$.
  \end{enumerate}

\end{lemma}

\begin{proof}
  Fix throughout the index $i$. We distinguish several steps.
  \paragraph{$\L1$-norm of $\Sigma_\rho w$.}
  With the change of variable $ \xi=X_i(0;t,x)$, for any
  $t \in [0,T]$,
  \begin{eqnarray*}
    &&\norma{(\Sigma_{\rho_o} w)_i(t)}_{\L1(\reali;\reali^n)}
    \\
    & =
    &\int_\reali
      \modulo{ (\rho_o)_i(X_i(0;t,x))
      \exp \left( -\int_{0}^t \partial_x w_i(s, X_i(s;t,x)) \, \d{s}\right)}
      \d{x}
    \\
    & =
    & \int_\reali \modulo{(\rho_o)_i(\xi)}
      \exp \left( -\int_{0}^t \partial_x w_i\left(s,X_i(s;t,x)\right) \, \d{s}\right)
      \exp \left( \int_0^t \partial_xw_i\left(s, X_i(s;0,\xi)\right)
      \d s \right) \d\xi
    \\
    & =
    & \int_\reali \modulo{(\rho_o(\xi))_i} \d\xi \,,
  \end{eqnarray*}
  proving~\eqref{eq:15}, thanks to~\eqref{eq:4} and~\eqref{eq:1}.

  \paragraph{$\tv$ estimate.}  If $\rho_o \in \W11(\reali;\reali^n)$
  then, for all $t \in [0,T]$, differentiating~\eqref{eq:4} \rev{and using
  the change of variable $ \xi = X_i(\tau;t,x)$,}
  \begin{eqnarray*}
    &&\norma{\partial_x (\Sigma_{\rho_o} w)_i (t)}_{\L1(\reali;\reali)}
    \\
    & =
    &
      \int_{\reali}
      \modulo{(\rho_o)_i'\left(X_i(0;t,x)\right)
      \partial_x X_i(0;t,x)
      \exp \left(-\int_0^t \partial_x w_i\left(\tau, X_i(\tau;t,x)\right) \d\tau
      \right)} \d{x}
    \\
    &
    & + \int_{\reali} \Big[\modulo{(\rho_o)_i\left(X_i(0;t,x)\right)}
      \exp \left( -\int_0^t \partial_x w_i\left(\tau, X_i(\tau;t,x)\right) \d\tau
      \right)
    \\
    &
    & \qquad \times
      \int_0^t \modulo{\partial^2_{xx} w_i\left(\tau, X_i(\tau;t,x)\right)
      \partial_xX_i(\tau;t,x)} \d\tau \Big] \d{x}
    \\
    & \leq
    & \exp \left( \norma{\partial_x w_i}_{\C0([0,t];\L\infty(\reali;\reali)}\, t \right) \tv\left((\rho_o)_i\right)
    \\
    &
    &+
      \norma{\partial_{xx}^2w_i}_{\C0([0,t];\L\infty(\reali))}
      \exp \left( \norma{\partial_x w_i}_{\C0([0,T];\L\infty(\reali;\reali)}t \right)
      \norma{(\rho_o)_i}_{\L1(\reali;\reali)}t \,,
  \end{eqnarray*}
  proving~\eqref{eq:tv_Sigma} when $\rho_o \in \W11(\reali;\reali^n)$.

  We now proceed to consider the case
  $\rho_o \in (\L1 \cap \BV) (\reali;\reali^n)$.  Thanks
  to~\cite[Theorem~3.9]{MR1857292}, there exists a sequence
  $\{ \rho_o^j\}_j \in \C\infty(\reali;\reali^n)$ such that
  $\rho_o^j \rightarrow \rho_o$ in $\L1(\reali;\reali^n)$ and
  $\lim_{j \rightarrow \infty} \norma{\partial_x \rho_o^j}_{\L1
    (\reali; \reali^n)} = \tv(\rho_o)$.

  Since, for $i = 1, \ldots, n$, \eqref{eq:15} ensures that
  \begin{equation*}
    \norma{(\Sigma_{\rho_o^j}w)_i(t) - (\Sigma_{\rho_o} w)_i(t)}_{\L1(\reali;\reali)}
    =
    \norma{(\rho_o^j)_i -(\rho_o)_i}_{\L1(\reali;\reali)} \,,
  \end{equation*}
  the lower semicontinuity of the total variation,
  see~\cite[Remark~3.5]{MR1857292}, ensures that
  \begin{eqnarray*}
    \tv\left((\Sigma_{\rho_o} w)_i(t)\right)
    &\leq
    & \liminf_{j \rightarrow \infty}
      \tv\left((\Sigma_{\rho_o^j}w)_i(t)\right)
    \\
    & \leq
    & \exp \left( \norma{\partial_x w_i}_{\C0([0,t];\L\infty(\reali;\reali)}
      \, t \right) \tv\left((\rho_o)_i\right)
    \\
    &+
    & \norma{\partial_{xx}^2w_i}_{\C0([0,t];\L\infty(\reali;\reali))}
      \exp \left( \norma{\partial_x w_i}_{\C0([0,t];\L\infty(\reali;\reali)} \, t
      \right)
      \norma{(\rho_o)_i}_{\L1(\reali;\reali)} t \,,
  \end{eqnarray*}
  completing the proof of~\eqref{eq:tv_Sigma}.

  \paragraph{$\Sigma_{\rho_o}$ is well defined.}  Given
  $w \in W_{\tilde Q}$ our aim is to prove the continuity in time of
  the function $\Sigma_{\rho_o} w$. So, fix $0 \leq t_1 < t_2 \leq T$,
  by~\eqref{eq:4} evaluate

\begin{align}
  \nonumber
  & \norma{ (\Sigma_{\rho_o} w)_i (t_1) - (\Sigma_{\rho_o} w)_i (t_2)}_{\L1(\reali;\reali)}
  \\\nonumber
  = & \int_\reali \Big| (\Sigma_{\rho_o} w)_i(t_1,x) - (\Sigma_{\rho_o} w)_i(t_1, X_i(t_1;t_2,x)) \exp\left(
      -\int_{t_1}^{t_2} \partial_xw_i(s;X_i(s;t_2,x))\d{s}
      \right) \Big| \dd x
  \\\nonumber
  \le &
        \int_\reali \Big| (\Sigma_{\rho_o} w)_i(t_1,x) - (\Sigma_{\rho_o} w)_i(t_1, X_i(t_1;t_2,x)) \Big| \dd x
  \\\nonumber
  & + \int_\reali \Big| (\Sigma_{\rho_o} w)_i(t_1, X_i(t_1;t_2,x)) \left( 1 - \exp\left(
    -\int_{t_1}^{t_2} \partial_xw_i(s;X_i(s;t_2,x))\dd s
    \right) \right) \Big| \dd x
  \\\label{eq:temp1}
  \leq & \tv\!\left((\Sigma_{\rho_o} w)_i(t_1)\right) \!\! \left( \!
         {\max} \! \left\{ \esssup_{x \in \reali} (x-X_i(t_1;t_2,x)),0\right\} {-} {\min} \! \left\{ \esssup_{x \in \reali} (x-X_i(t_1;t_2,x)),0\right\} \! \right)
  \\\label{eq:temp2}
  & + \int_\reali \Big| (\Sigma_{\rho_o} w)_i(t_1;X_i(t_1;t_2,x)) \int_{t_1}^{t_2} \partial_x w_i(s;X_i(s;t_2,x))\dd s \Big| \,\dd x,
\end{align}
where we used \Cref{lem:TV}.

Consider~\eqref{eq:temp1}. Using the estimate~\eqref{eq:tv_Sigma} and
the boundedness of $w_i$ we obtain
\begin{equation*}
  \tv((\Sigma_{\rho_o} w)_i(t_1)) \leq e^{{\tilde Q}t_1} \tv((\rho_o)_i) +{\tilde Q}t_1e^{{\tilde Q}t_1} \norma{(\rho_o)_i }_{\L1(\reali;\reali)}
\end{equation*}
and
\begin{eqnarray*}
  &
  & \max \left\{ \esssup\nolimits_{x \in \reali} (x-X_i(t_1;t_2,x)),0\right\} - \min\left\{ \esssup\nolimits_{x \in \reali} (x-X_i(t_1;t_2,x)),0\right\}
  \\
  & \leq
  & 2 \esssup_{x \in \reali} \modulo{x - X_i(t_1;t_2,x)}
  \\
  & \leq
  & 2 \norma{w_i}_{\C0([0,T];\L\infty(\reali;\reali))} (t_2-t_1)
  \\
  & \leq
  & 2 {\tilde Q}(t_2-t_1)
\end{eqnarray*}
so that
\begin{displaymath}
  [\eqref{eq:temp1}]
  \leq
  \left( e^{{\tilde Q}t_1} \tv((\rho_o)_i) +{\tilde Q}t_1e^{{\tilde Q}t_1} \norma{(\rho_o)_i }_{\L1(\reali;\reali)} \right)
  2{\tilde Q}\, (t_2-t_1) \,.
\end{displaymath}

Consider now~\eqref{eq:temp2}.  By the change of variable
$\xi= X_i(t_1;t_2,x)$, \eqref{eq:15} and~\eqref{eq:5}
\begin{eqnarray*}
  {[\eqref{eq:temp2}]}
  & \le
  & {\tilde Q}(t_2 -t_1) \int_\reali \Big| (\Sigma_{\rho_o} w)_i(t_1;X_i(t_1;t_2,x)) \Big| \,\dd x
  \\
  & \le
  & {\tilde Q}(t_2 -t_1) \int_\reali \big| \Sigma_{\rho_o} w(t_1, \xi)_i \big| \exp\left( \int_{t_1}^{t_2} \partial_xw_i(s, X_i(s;t_1, \xi)) \,\d{s} \right) \, \dd\xi
  \\
  & \leq
  & {\tilde Q}(t_2-t_1) e^{{\tilde Q} (t_2-t_1)}
    \norma{(\rho_o)_i}_{\L1(\reali;\reali)}.
\end{eqnarray*}
Collecting together the estimates above leads to the following bound:
\begin{eqnarray}
  \label{eq:25}
  &
  &\norma{(\Sigma_{\rho_o} w)_i (t_1) - (\Sigma_{\rho_o}w)_i(t_2)}_{\L1(\reali;\reali)}
  \\
  \nonumber
  & \leq
  & {\tilde Q} \;
    \left[
    2 \left(\tv((\rho_o)_i) +{\tilde Q} T \norma{ (\rho_o)_i}_{\L1(\reali;\reali)}\right)
    +
    \norma{(\rho_o)_i}_{\L1(\reali;\reali)}
    \right]     e^{{\tilde Q} T}
    (t_2-t_1)
\end{eqnarray}
concluding the proofs of~\ref{item:9} and of the well-posedness of
$\Sigma_{\rho_o}$.

\paragraph{The map $\Sigma_{\rho_o}$ in~\eqref{def:Sigma} is Lipschitz
  continuous.}  Consider $t\in [0,T]$, $w,z \in W_{\tilde Q}$ and
using~\eqref{eq:4}, estimate
\begin{eqnarray}
  \label{eq:Lip_Sigma}
  \!\!\!\!\!\!
  \norma{(\Sigma_{\rho_o} w)_i (t)- (\Sigma_{\rho_o} z)_i (t) }_{\L1(\reali;\reali)} \!\!\!
  & \!\!\!\!\!\!\!\!\!=\!\!\!\!\!\!\!\!\!
  & \!\!\! \int_\reali \Big| \rho_o(X_i(0;t,x))_i \exp \left(\!{-}\!
    \int_0^{t} \partial_xw_i(\tau,X_i(\tau;t,x)) \d\tau
    \right)
  \\
  \nonumber
  &
  & -
    \rho_o(Y_i(0;t,x))_i \exp \left( -
    \int_0^{t} \partial_xz_i(\tau,Y_i(\tau;t,x)) \dd \tau
    \right) \Big| \, \dd x
\end{eqnarray}
where $X_i(t;t_o,x_o)$ and $Y_i(t;t_o,x_o)$ are the solutions to
\begin{equation}
  \begin{cases}
    \dot{X_i} = w_i (t, x)\\
    X_i(t_o)=x_o
  \end{cases}
  \quad \mbox{ and } \quad
  \begin{cases}
    \dot{Y_i} = z_i(\tau,y)\\
    Y_i(t_o)=x_o.
  \end{cases}
\end{equation}
Adding and subtracting the term
$\rho_o(X_i(0;t,x))_i\exp\left(-\int_0^{t}
  \partial_xz_i(\tau,Y_i(\tau;t,x)) \dd \tau \right)$
to~\eqref{eq:Lip_Sigma} and thanks to~\eqref{eq:5},
Proposition~\ref{prop 2} and Lemma~\ref{lem:TV}, one obtains
\begin{align*}
  &\quad\norma{(\Sigma_{\rho_o}w)_i(t) - (\Sigma_{\rho_o}z)_i(t) }_{\L1(\reali;\reali)}
  \\
  & \leq
    \int_\reali |\rho_o(X_i(0;t,x))_i - \rho_o(Y_i(0;t,x))_i| \dd x
    e^{{\tilde Q}\,t}
  \\
  &\quad
    +\int_\reali |\rho_o(X_i(0;t,x))_i| \Big| \int_0^t \partial_xw_i(\tau, X_i(\tau;t,x)) - \partial_xz_i(\tau, Y_i(\tau;t,x)) \dd\tau \Big| \dd x
    e^{{\tilde Q}\,t}
  \\
  & \leq
    2 \, \tv((\rho_o)_i) \,
    \norma{X_i(0;t,\cdot) - Y_i(0;t,\cdot)}_{\L\infty(\reali;\reali)} \,
    e^{{\tilde Q}\,t}
  \\
  & \quad {+} \!\!\! \int_\reali \modulo{\rho_o(X_i(0;t,x))_i}
    \left(\!
    {\tilde Q} \!\! \int_0^t \!\! \modulo{X_i(\tau;t,x) {-} Y_i(\tau;t,x)} \dd\tau {+}
    \norma{\partial_xw_i {-} \partial_xz_i}_{\C0([0,T];\L\infty(\reali;\reali))}t
    \!\right) \! \d{x} e^{{\tilde Q}\,t}
  \\
  & \leq
    2\, \tv((\rho_o)_i) \,
    \norma{w_i - z_i}_{\C0([0,T]; \L\infty(\reali;\reali))} \, t \, e^{2\,{\tilde Q}\,t}
  \\
  & \quad
    + \int_\reali \modulo{\rho_o(X_i(0;t,x))_i} \left[
    {\tilde Q} \int_0^t \norma{w_i -z_i}_{\C0([0,T]; \L\infty(\reali;\reali))} \, |t-\tau|
    \, e^{{\tilde Q}\,(t-\tau)} \, \dd \tau
    \right ] \, \dd x \, e^{{\tilde Q}\,t}
  \\
  & \quad
    + \int_\reali \modulo{\rho_o(X_i(0;t,x))_i}
    \left[ t
    \norma{\partial_x w_i - \partial_x z_i}_{\C0([0,T]; \L\infty(\reali;\reali))}
    \right ] \, \dd x \, e^{{\tilde Q}\,t}
  \\
  & \leq
    2\, \tv((\rho_o)_i) \, \norma{w_i - z_i}_{\C0([0,T]; \L\infty(\reali;\reali))}
    \, t \, e^{2\,{\tilde Q}\,t}
  \\
  & \quad
    + {\tilde Q} \, \norma{(\rho_o)_i}_{\L1(\reali;\reali)}\, \norma{w_i - z_i}_{\C0([0,T]; \L\infty(\reali;\reali))}
    \, \frac{t^2}{2}
    \, e^{2\, {\tilde Q} \, t}
  \\
  & \quad
    +
    \norma{(\rho_o)_i}_{\L1(\reali;\reali)}\, \norma{\partial_x w_i - \partial_x z_i}_{\C0([0,T]; \L\infty(\reali;\reali))}\, t \, e^{2\,{\tilde Q}\,t}
  \\
  & \leq
    \left(
    2\tv((\rho_o)_i)
    +
    \frac{1}{2} \, {\tilde Q} \, t \norma{(\rho_o)_i}_{\L1(\reali;\reali)}
    + \norma{(\rho_o)_i}_{\L1(\reali;\reali)} \right)
    \norma{w_i -z_i}_{\C0([0,T]; \W1\infty(\reali;\reali))} t \, e^{2\,{\tilde Q}\,t}.
\end{align*}
Passing to the supremum over $t \in [0,T]$,
\begin{equation*}
  \norma{ \Sigma_{\rho_o} w - \Sigma_{\rho_o} z}_{\C0([0,T];\L1(\reali;\reali^n))} \leq
  \Lip(\Sigma_{\rho_o}) \, \norma{w - z}_{\C0([0,T];\W1\infty(\reali;\reali^n))}
\end{equation*}
where
\begin{equation}
  \label{eq:17}
  \Lip(\Sigma_{\rho_o})
  \coloneqq
  \left(
    2\tv(\rho_o)
    +
    \frac{1}{2} \, {\tilde Q} \, T \norma{\rho_o}_{\L1(\reali;\reali)}
    + \norma{\rho_o}_{\L1(\reali;\reali)} \right)   T\, \exp(2\,{\tilde Q}\,T).
\end{equation}
The proof of~\ref{item:3} is completed.

Finally, the positivity~\ref{item:13} immediately follows
from~\eqref{eq:4}.
\end{proof}

\begin{proposition}
  \label{prop:4}
  Let $\tilde Q >0$. For a fixed $w \in W_{\tilde Q}$ as defined
  in~\eqref{eq:5}, the map
  \begin{equation}
    \begin{array}{ccc}
      (\L1 \cap \BV) (\reali; \reali^n)
      & \to
      & \C0([0,T];\L1(\reali;\reali^n))
      \\
      \rho_o
      & \mapsto
      & \Sigma_{\rho_o} w
    \end{array}
  \end{equation}
  with $\Sigma$ as in~\eqref{def:Sigma}, is Lipschitz continuous.
\end{proposition}

\begin{proof}
  Choose $\rho_o$ and $\sigma_o$ in
  $(\L1 \cap \BV) (\reali; \reali^n)$. Fix $w \in W_{\tilde Q}$, $t \in [0,T]$
  and $i = 1, \ldots, n$.  Then the result is a direct consequence of
  the linearity of~\eqref{eq:14} and of the equality~\eqref{eq:15}:
  \begin{displaymath}
    \norma{\Sigma_{\rho_o}w(t) - \Sigma_{\sigma_o}w(t)}_{\L1(\reali;\reali^n)}
    =
    \norma{\Sigma_{\rho_o-\sigma_o} w(t)}_{\L1(\reali;\reali^n)}
    =
    \norma{\rho_o - \sigma_o}_{\L1(\reali;\reali^n)} \,.
  \end{displaymath}
\end{proof}

\begin{lemma}\label{lem:mathcalt}
  Let~\ref{item:1} and~\ref{item:2} hold. Fix $M>0$ and
  $\rho_o \in (\L1 \cap \BV) (\reali; \reali^n)$ such that
  \begin{equation}
    \label{stima_M}
    M  = \norma{\rho_o}_{\L1(\reali;\reali^n)} \,.
  \end{equation}
  Define $Q \coloneqq \sup_{t \in [0,T]} Q_{v (t)}$ as
  in~\eqref{stima_Q} and $X_M$ as
  in~\eqref{eq:10}--\eqref{eq:2}. Then, the map
  \begin{equation}\label{mathcal_t}
    \begin{array}{cccc}
      \mathcal{T}_{v, \eta, \rho_o}:
      &  X_M
      & \rightarrow
      & X_M
      \\
      & \rho
      & \mapsto
      & \left(\Sigma_{\rho_o } \circ \Pi_{v, \eta}\right) (\rho)
    \end{array}
  \end{equation}
  is well defined and Lipschitz continuous. For each $\rho \in X_M$
  the function $\mathcal{T}_{v, \eta, \rho_o} \rho$ is locally Lipschitz
  continuous in time and with total variation in space bounded by
  \begin{equation}
    \label{eq:19}
    \tv\left(\mathcal{T}_{v, \eta, \rho_o} \rho(t)\right)
    \leq
    \exp(Qt) \left( \tv(\rho_o) +
      Q \, \norma{\rho_o}_{\L1(\reali;\reali^n)} \, t\right)
    \quad \forall \, t \in [0,T].
  \end{equation}
  Furthermore, if $T$ is sufficiently
  small, then $\mathcal{T}_{v, \eta, \rho_o}$ is also a contraction.
\end{lemma}

\begin{proof}
  Throughout this proof, we keep $v$, $\eta$ and $\rho_o$ fixed, hence
  we omit them.

  \paragraph{The map $\mathcal{T}$ is well defined and Lipschitz
    continuous.} Thanks to the property of the maps $\Pi$ and
  $\Sigma$, it is immediate to verify that if $\rho \in X$ then
  $\mathcal{T} \rho \in \C0([0,T];\L1(\reali;\reali^n))$.
  Furthermore, owing to~\eqref{eq:15},
  $\norma{\mathcal{T}\rho }_{\C0([0,T];\L1(\reali;\reali^n))} =
  \norma{\rho_o}_{\L1(\reali,\reali^n)} = M$.

  The Lipschitz regularity follows, since $\mathcal{T}$ is the
  composition of Lipschitz continuous maps, by Lemma~\ref{lem:Sigma}
  and Lemma~\ref{lem:Pi}.

    \paragraph{Properties of the function $\mathcal{T}\rho.$}
    Given $\rho \in X_M$, the local $\L1$-Lipschitz continuity in time
    of $\mathcal{T}\rho$ directly descends from~\ref{item:9} in
    \Cref{lem:Sigma}.  Furthermore, for each $t \in [0,T]$, thanks
    to~\ref{eq:tv_Sigma}, one obtains that
    \begin{displaymath}
      \tv(\mathcal{T}\rho(t))
      \leq
      \exp(Q \, t) \left( \tv(\rho_o) +
        Q \, \norma{\rho_o}_{\L1(\reali;\reali^n)} \, t\right)
    \end{displaymath}
    proving~\eqref{eq:19}.

    \paragraph{The map $\mathcal{T}$ is a contraction for sufficiently
      small times.} Since we have
    $\Lip(\mathcal{T}) = \Lip(\Sigma) \Lip(\Pi)$, we are lead to prove
    that $\Lip(\Sigma) \Lip(\Pi) < 1$.

    Choose $T < 2$, which implies that
    \begin{align*}
      Q
      & = \sup_{t \in [0, T]} Q_{v(t)} \leq C \left(\norma{\eta}_{\W{2}{\infty}(\reali; \reali^{n \times m})}, M, \sup_{t \in [0, 2]}  \norma{v(t)}_{\mathcal{V}^n}\right) =: \bar Q.
    \end{align*}
    Moreover, by~\eqref{eq:16},
    \begin{equation*}
      \Lip (\Pi) \le C\left(\norma{\eta}_{\W{2}{\infty}(\reali; \reali^{n \times m})}, M, \sup_{t \in [0, 2]}  \norma{v(t)}_{\mathcal{V}^n}\right) =: \overline{\Lip \Pi}.
    \end{equation*}

    Hence, by~\eqref{eq:17} the additional conditions
    \begin{eqnarray*}
      T
      & <
      & 1/ (2\,\bar Q)
      \\
      T
      & <
      & \dfrac{1}{2}
        \left(
        e
        \left(
        2 \, \tv((\rho_o)_i)
        +
        (\bar Q+1) \, M
        \right)
        \overline{\Lip(\Pi)}
        \right)^{-1}
    \end{eqnarray*}
    ensure that $\Lip(\mathcal{T}) \leq 1/2$.
  \end{proof}

\begin{proposition}\label{prop:dipendenza_lip}
  Fix $\rho \in X_M$ as defined in~\eqref{eq:10}.  Define the map
  \begin{equation}\label{eq:18}
    \begin{array}{ccc}
      \C0([0,T];\mathcal{V}^n) \times \W2\infty(\reali;\reali^{n\times n})\times (\L1 \cap \BV)(\reali;\reali^n)
      & \rightarrow
      & \C0\left([0,T];\L1(\reali;\reali^n)\right)
      \\
      v, \eta, \rho_o
      & \mapsto
      & \mathcal{T}_{v, \eta, \rho_o}\rho.
    \end{array}
  \end{equation}
  Then
  \begin{enumerate}[label=\textbf{(\arabic*)}]
  \item For all $\eta \in \W2\infty(\reali;\reali^{n \times n })$,
    $\rho_o \in (\L1 \cap \BV)(\reali;\reali^n)$, the map
    $v \mapsto \mathcal{T}_{v, \eta, \rho_o}\rho$ is Lipschitz
    continuous.
  \item For all $v \in \C0([0,T];\mathcal{V}^n)$,
    $\rho_o \in (\L1 \cap \BV)(\reali;\reali^n)$, the map
    $ \eta \mapsto \mathcal{T}_{v, \eta, \rho_o}\rho$ is locally
    Lipschitz continuous.
  \item For all $ v \in \C0([0,T];\mathcal{V}^n)$,
    $\eta \in \W2\infty(\reali;\reali^{n \times n})$ the map
    $\rho_o \mapsto \mathcal{T}_{v, \eta, \rho_o}\rho$ is Lipschitz
    continuous.
  \end{enumerate}
\end{proposition}

\noindent The proof follows directly from Proposition~\ref{prop:5}
and Proposition~\ref{prop:4}.

\begin{proofof}{Theorem~\ref{thm:main}}

  \paragraph{Proof of~\ref{item:4}:}
  Note that Definition~\ref{def:solution} implies that
  solving~\eqref{eq:9} is equivalent to finding a fixed point of the
  map $\mathcal{T}_{v, \eta, \rho_o}$ defined
  in~\ref{mathcal_t}. Lemma~\ref{lem:mathcalt} ensures that for a
  $T^*>0$ the map $\mathcal{T}_{v, \eta, \rho_o}$ admits a unique
  fixed point.

  Fix an arbitrary $T > 0$. Recall the quantities $M$ defined
  in~\eqref{stima_M}, $Q \coloneqq \sup_{t \in [0,T]} Q_{v (t)}$ as
  in~\eqref{stima_Q} and by~\eqref{eq:16}
  \begin{equation*}
    \Lip (\Pi_{v,\eta})
    \le
    C\left(
      \norma{\eta}_{\W{2}{\infty}(\reali; \reali^{n \times m})},
      M,
      \sup_{t \in [0, T]}  \norma{v(t)}_{\mathcal{V}^n}
    \right)
    =: \overline{\Lip(\Pi_{v, \eta})}.
  \end{equation*}
  Introduce
  \begin{eqnarray*}
    K
    & \coloneqq
    &  \exp (Q \, T) \left( \tv(\rho_o) + Q\,M\, T\right)
    \\
    \Delta T
    & \coloneqq
    & \min\left\{T, \frac{1}{2\, e \left(2\,K + Q\,M + M\right) \overline{\Lip(\Pi_{v,\eta})}}\right\} \,.
  \end{eqnarray*}
  Then, Lemma~\ref{lem:mathcalt} ensures that there exists
  $\rho \in \C0([0,\Delta T]; \L1(\reali;\reali^n))$
  solving~\eqref{eq:9} on $[0,\Delta T]$. Observe that
  $\norma{\rho(\Delta T)}_{\L1(\reali;\reali^n)}=M$ by~\eqref{eq:15}
  and $\tv\left(\rho(\Delta T)\right) \leq K$ by~\eqref{eq:19}. We can
  iterate the application of Lemma~\ref{lem:mathcalt} obtaining a
  solution to~\eqref{eq:9} on $[0,T]$, since
  $\norma{\rho(k\, \Delta T)}_{\L1(\reali;\reali^n)} = M$ and
  $\tv\left(\rho(k\, \Delta T)\right) \leq K$ for $k=1,2,\ldots$. By
  the arbitrariness of $T$, \ref{item:4} follows.

  Define $\mathcal{P}_{0,t} \rho_o = \rho (t)$.

  \paragraph{Proof of~\ref{item:16}:} This property directly follows
  from the construction, for instance from~\eqref{eq:4}.

  \paragraph{Proof of~\ref{item:5}:}
  Call
  $M \coloneqq \max\{\norma{\rev{\rho_o}}_{\L1(\reali;\reali^n)},
  \norma{\rev{\hat \rho_o}}_{\L1(\reali;\reali^n)}\}$. For
  $t \in \reali_+$ denote $\rho'(t) = \mathcal{P}_{0,t}\rev{\rho_o}$,
  $\rho''(t) = \mathcal{P}_{0,t}\rev{\hat \rho_o}$ and introduce the
  corresponding characteristics $X'_i$ and $X''_i$ as
  in~\eqref{eq:24}. Then, by~\eqref{mathcal_t}
  $ \rho' = \mathcal{T}_{v , \eta, \rev{\rho_o}}\rho'$ and
  $\rho'' = \mathcal{T}_{v , \eta, \rev{\hat \rho_o}}\rho''$, so that
  for $i=1,\ldots, n$
  \begin{align*}
    &
      \norma{(\rho')_i(t)- (\rho'')_i(t)}_{\L1(\reali;\reali)}
    \\
    & \leq
      \int_\reali \modulo{(\rev{\rho_o})_i(X'_i(0;t,x)) - (\rev{\hat \rho_o})_i(X'_i(0;t,x))} \exp \left( -\int_0^t \partial_x (\Pi_{v, \eta} \rho')_i(s, X'_i(s;t,x)) \dd s \right) \dd x
    \\
    & +
      \int_\reali \modulo{(\rev{\hat \rho_o})_i(X'_i(0;t,x)) - (\rev{\hat \rho_o})_i({X}_i''(0;t,x))} \exp \left( -\int_0^t \partial_x (\Pi_{v, \eta} \rho')_i(s, X'_i(s;t,x)) \dd s \right) \dd x
    \\
    & +
      e^{Qt}\int_\reali \modulo{(\rev{\hat \rho_o})_i({X}_i''(0;t,x))}  \modulo{ \int_0^t\partial_x (\Pi_{v, \eta} \rho')_i(s, X'_i(s;t,x))- \partial_x (\Pi_{v, \eta} {\rho''})_i(s, {X}_i''(s;t,x)) \dd s} \dd x
    \\
    &\leq
      \norma{(\rev{\rho_o})_i - (\rev{\hat \rho_o})_i}_{\L1 (\reali;\reali)} + A_1 + A_2 + A_3 \,,
  \end{align*}
  where
  $Q = C(\norma{\eta}_{\W2\infty(\reali;\reali^{n\times n})}, \tilde
  M, \norma{v}_{\C0([0,t];\mathcal{V}^n)}$) as in \Cref{lem:mathcalt},
  and
  \begin{align*}
    A_1
    & =
      2 \, e^{Qt} \tv\left((\rev{\hat \rho_o})_i\right) \, \norma{X_i'(0;t,\cdot) - {X}_i''(0;t, \cdot)}_{\L\infty(\reali;\reali)}
    \\
    A_2
    & =
      e^{Qt} \!\!\int_\reali \modulo{(\rev{\hat \rho_o})_i({X}_i''(0;t,x))}  \modulo{ \int_0^t\partial_x (\Pi_{v, \eta} \rho')_i(s, X_i'(s;t,x)) {-} \partial_x (\Pi_{v, \eta} {\rho''})_i(s, X_i'(s;t,x)) \dd s} \dd x
    \\
    A_3
    & =
      e^{Qt} \!\!\int_\reali \modulo{(\rev{\hat \rho_o})_i({X}_i''(0;t,x))}  \modulo{ \int_0^t\partial_x (\Pi_{v, \eta} {\rho''})_i(s, X_i'(s;t,x)) {-} \partial_x (\Pi_{v, \eta} {\rho''})_i(s, {X}_i''(s;t,x)) \dd s} \dd x.
  \end{align*}
  We estimate the latter terms separately:
  \begin{flalign*}
    & \norma{X_i'(0;t,\cdot) - {X}_i''(0;t,
      \cdot)}_{\L\infty(\reali;\reali)}
    \\
    \leq & \norma{\partial_x\left( \Pi_{v,\eta} \rho'\right)_i -
      \partial_x\left( \Pi_{v,\eta} \rho''\right)_i}_{\L1
      ([0,t];\L\infty (\reali; \reali))} &
    [\mbox{By~\eqref{eq:14}--\eqref{eq:24}}]
    \\
    \leq & \Lip (\Pi_{v,\eta}) \, \norma{\rho'-\rho''}_{\L1
      ([0,t]\times \reali; \reali^n)} & [\mbox{By
      Lemma~\eqref{lem:Pi}]}
  \end{flalign*}
  so that
  \begin{displaymath}
    A_1
    \leq
    2 \, \Lip (\Pi_{v,\eta}) \, \tv\left((\rev{\hat \rho_o})_i\right)
    \, e^{Q\,t}
    \, \norma{\rho'-\rho''}_{\L1 ([0,t]\times\reali; \reali^n)} \,.
  \end{displaymath}
  Moreover,
  \begin{flalign*}
    & \modulo{ \int_0^t\partial_x (\Pi_{v, \eta} \rho')_i(s,
      X_i'(s;t,x))- \partial_x (\Pi_{v, \eta} {\rho''})_i(s,
      X_i'(s;t,x)) \dd s}
    \\
    \leq & \int_0^t \norma{\partial_x (\Pi_{v, \eta} \rho')_i(s,
      X_i'(s;t,x))- \partial_x (\Pi_{v, \eta} {\rho}'')_i(s,
      X_i'(s;t,x))}_{\L\infty (\reali; \reali^n)} \dd s
    \\
    \leq & \int_0^t\Lip (\Pi_{v,\eta}) \, \norma{\rho' (s) -
      \rho''(s)}_{\L1 (\reali; \reali^n)} \d{s} & [\mbox{By
      Lemma~\eqref{lem:Pi}]}
  \end{flalign*}
  implying that
  \begin{displaymath}
    A_2
    \leq
    \Lip (\Pi_{v,\eta}) \, M \, e^{2\, Q \, t}
    \, \norma{\rho' - \rho''}_{\L1 ([0,t]\times\reali; \reali^n)}  \,.
  \end{displaymath}
  Finally,
  \begin{flalign*}
    & \modulo{ \int_0^t\partial_x (\Pi_{v, \eta} {\rho''})_i(s,
      X_i'(s;t,x))- \partial_x (\Pi_{v, \eta} {\rho''})_i(s,
      {X}_i''(s;t,x)) \dd s}
    \\
    \leq & \int_0^t \modulo{\partial_x (\Pi_{v, \eta} {\rho''})_i(s,
      X'_i(s;t,x))- \partial_x (\Pi_{v, \eta} {\rho''})_i(s,
      {X}_i''(s;t,x))} \dd s
    \\
    \leq & \int_0^t \norma{\Pi_{v,\eta}\rho''}_{\W2\infty (\reali;
      \reali^n)} \, \modulo{X_i' (s;t,x) - X_i'' (s;t,x)}\d{s}
    \\
    \leq & Q \, \Lip (\Pi_{v,\eta}) \, t \, \norma{\rho'-\rho''}_{\L1
      ([0,t]\times \reali; \reali^n)}
  \end{flalign*}
  so that
  \begin{displaymath}
    A_3
    \leq Q \, M \, \Lip (\Pi_{v,\eta}) \, t \, e^{2\,Q\,t} \, \norma{\rho'-\rho''}_{\L1 ([0,t]\times \reali; \reali^n)} \,.
  \end{displaymath}
  We thus obtain
  \begin{eqnarray*}
    &
    & \norma{\rho'(t)- {\rho}''(t)}_{\L1(\reali;\reali^n)}
    \\
    & \leq
    & \norma{\rev{\rho_o} - \rev{\hat \rho_o}}_{\L1 (\reali;\reali^n)}
    \\
    &
    & +\left(2  \, \tv\left(\rev{\hat \rho_o}\right)
      + M
      + Q \, M \, t \right)  \, \Lip (\Pi_{v,\eta}) \, e^{2\, Q \, t}
      \norma{\rho'-\rho''}_{\L1 ([0,t]\times \reali; \reali^n)} \,.
  \end{eqnarray*}
  An application of Gronwall Lemma completes the proof
  of~\ref{item:5}.

  \paragraph{Proof of~\ref{item:8}:}
  Choose $t',t'' \in [0,T]$. Then, calling
  $w (t,x) = \Pi_{v,\eta}\rho (t,x)$,
  \begin{flalign*}
    & \norma{\rho (t'')-\rho (t')}_{\L1 (\reali; \reali^n)}
    \\
    = & \norma{(\Sigma_{\rho_o} w) (t'') - (\Sigma_{\rho_o} w)
      (t')}_{\L1 (\reali; \reali^n)} & [\mbox{Since }
    \mathcal{T}_{v_\eta,\rho_o} \rho = \rho]
    \\
    \leq & C\left(\norma{\eta}_{\W2\infty (\reali;\reali^{n\times
          n})}, \norma{v}_{\C0([0,T];\mathcal{V}^n)},
      \norma{\rho_o}_{\L1 (\reali;\reali^n)}, \tv (\rho_o) , T\right)
    \; \modulo{t'' - t'} & [\mbox{By~\eqref{eq:25} in
      Lemma~\ref{lem:Sigma}}]
  \end{flalign*}
  proving~\ref{item:8}.

  \paragraph{Proof of~\ref{item:6}:} For any $t \in \reali_+$,
  by~\cite[Theorem~2.9]{MR1816648}, we have:
  \begin{displaymath}
    \norma{\rev{\mathcal{P}}_{0,t}\rho_o - \rev{\hat{\mathcal{P}}}_{0,t}\rho_o}_{\L1 (\reali;\reali^n)}
    \leq
    \Lip (\rev{\mathcal{P}})
    \int_0^t \liminf_{h\to 0^+} \dfrac{\norma{\rev{\mathcal{P}}_{\tau,h} \rev{\hat{\mathcal{P}}}_{0,\tau} \rho_o - \rev{\hat{\mathcal{P}}}_{\tau,h} \rev{\hat{\mathcal{P}}}_{0,\tau} \rho_o}_{\L1 (\reali;\reali^n)}}{h} \d\tau \,
  \end{displaymath}
  where
  \begin{equation*}
    \Lip(\rev{\mathcal{P}})
    \coloneqq
    C\left(\norma{\eta}_{\W2\infty(\reali;\reali^{n \times n})},
      \norma{\rev{v}}_{\C0([0,t];\mathcal{V}^n)},
      \norma{\rho_o}_{\L1(\reali;\reali^n)}, \tv(\rho_o), t\right)
  \end{equation*}
  is as in~\ref{item:5}.  Call
  $\rho = \rev{\hat{\mathcal{P}}}_{0,\tau} \rho_o$. Recall that for
  $h$ small, owing to Proposition~\ref{prop:dipendenza_lip},
  \begin{displaymath}
    \norma{\rev{\mathcal{P}}_{\tau,h}\rho - \rev{\hat{\mathcal{P}}}_{\tau,h}\rho}_{\L1 (\reali; \reali^n)}
    \leq \frac{B_1 h}{1-B_2h} \; \norma{\rev{v}-\rev{\hat v}}_{\C0([0, h]; \mathcal{V}^n)} \;
  \end{displaymath}
  with
  \begin{eqnarray*}
    B_1
    & \coloneqq
    & \left(2\tv(\rho_o)+\left(\frac{1}{2}Qh + 1\right) \norma{\rho_o}_{\L1(\reali; \reali^n)}\right) \;
      C\left(\norma{\rho_o}_{\L1(\reali;\reali^n)}, \norma{\eta}_{\W2\infty(\reali;\reali^{n \times n})}\right)
    \\
    B_2
    & \coloneqq
    & \left( 2\tv(\rho_o)+\left(\frac{1}{2}Q h + 1 \right) \norma{\rho_o}_{\L1(\reali; \reali^n)}\right)\Lip(\Pi_{\rev{\hat v},\eta})
  \end{eqnarray*}
  and the constant
  $C(\norma{\rho_o}_{\L1(\reali;\reali^n)},
  \norma{\eta}_{\W2\infty(\reali;\reali^{n \times n})})$ is given
  by~\ref{item:10} in Proposition~\ref{prop:5}.  So,
  \begin{eqnarray*}
    \!\!\!
    &
    & \norma{\rev{\mathcal{P}}_{0,t}\rho_o - \rev{\hat{\mathcal{P}}}_{0,t}\rho_o}_{\L1 (\reali;\reali^n)}
    \\
    \!\!\!
    & \leq
    & \left( 2\tv(\rho_o) + \norma{\rho_o}_{\L1(\reali; \reali^n)} \right)
      C \left(\norma{\rho_o}_{\L1(\reali; \reali^n)},
      \norma{\eta}_{\W2\infty(\reali;\reali^{n \times n})} \right) \,
      \Lip(\rev{\mathcal{P}})\, t\,  \norma{\rev{v}-\rev{\hat v}}_{\C0([0,t];\mathcal{V}^n)}.
  \end{eqnarray*}
  completing the proof of~\ref{item:6}.

\paragraph{Proof of~\ref{item:7}:}
For any $ t \in \reali_+$, using the same technique as in the proof
of~\ref{item:6}, one gets
\begin{align*}
  & \norma{\rev{\mathcal{P}}_{0,t}\rho_o - \rev{\hat{\mathcal{P}}}_{0,t}\rho_o}_{\L1 (\reali;\reali^n)}
  \\
  \leq
  & C\left( \norma{\rev{\eta}}_{\W2\infty(\reali;\reali^{n \times n})}, \norma{\rev{\hat \eta}}_{\W2\infty(\reali;\reali^{n \times n})}, \norma{\rho_o}_{\L1(\reali,\reali^n)}\right)\\
  & \times \left( 2\tv(\rho_o) + \norma{\rho_o}_{\L1(\reali,\reali^n)} \right) \norma{\rho_o}_{\L1(\reali,\reali^n)}\norma{v}_{\C0([0,t];\mathcal{V}^n)}
    \Lip(\rev{\mathcal{P}})\, t\,  \norma{\rev{\eta}-\rev{\hat \eta}}_{\W2\infty(\reali;\reali^{n \times n})}
\end{align*}
thanks to Proposition~\ref{prop:dipendenza_lip}.

The total variation estimate~\ref{item:15} directly follows
from~\eqref{eq:19} while the positivity~\ref{item:14} is
immediate. Hence, the proof is completed.
\end{proofof}

\section*{Acknowledgments}
RMC and MG were partly supported by the GNAMPA~2023 project
\emph{Analytical Techniques for Biological Models, Fluid dynamics and
  Vehicular Traffic}. They also acknowledge the PRIN 2022 project
\emph{Modeling, Control and Games through Partial Differential
  Equations} (CUP~D53D23005620006), funded by the European Union -
Next
Generation EU.\\
MG was partly supported by the European Union-NextGeneration EU
(National Sustainable Mobility Center CN00000023, Italian Ministry of
University and Research Decree n.~1033-17/06/2022, Spoke 8) and by the
Project funded under the National Recovery and Resilience Plan (NRRP)
of Italian Ministry of University and Research funded by the European
Union-NextGeneration EU. Award Number: ECS\_00000037, CUP:
H43C22000510001, MUSA-Multilayered Urban Sustainability Action.\\
CN was supported by the PNRR project \emph{Ricerca DM
  118/2023}, CUP~F13C23000360007.

\section*{Conflict  of Interest}
The authors declare no conflicts of interest in this paper.

{\small

  \bibliographystyle{abbrv}

  \bibliography{Nonlocal_time_dependent_speed}

}

\end{document}